\documentclass[a4paper,reqno,11pt]{amsart}
\usepackage[T1]{fontenc}
\usepackage[utf8]{inputenc}
\usepackage{amssymb, amsmath, amsthm, graphicx, color}
\usepackage[inline]{enumitem}
\usepackage[dvipsnames,svgnames,table]{xcolor}
\usepackage[foot]{amsaddr}
\usepackage{comment}

\usepackage{booktabs}
\usepackage{multirow}

\usepackage[longnamesfirst,numbers,sort&compress]{natbib}

\providecommand{\noopsort}[1]{}

\usepackage{comment}

\makeatletter
\def\thm@space@setup{
\thm@preskip=4mm
\thm@postskip=0mm
}
\makeatother

\usepackage{hyperref}
\hypersetup{
    pdftitle={Centered colorings in minor-closed graph classes},
    pdfauthor={Jędrzej Hodor, Hoang La, Piotr Micek, Clément Rambaud},
    colorlinks,
    linkcolor={RoyalBlue},
    citecolor={RubineRed},
    urlcolor={blue!80!black}
}

\usepackage[capitalise, noabbrev]{cleveref}

\usepackage{mathtools}

\usepackage{thmtools, thm-restate}

\setenumerate{label=\textup{(\roman*)}, noitemsep, topsep=3pt-\parskip, 
labelindent=.2em, leftmargin=*, widest=iii,}
\setitemize{noitemsep, topsep=-\parskip, labelindent=.2em, leftmargin=*, widest=iii,}

\newenvironment{enumerateOurAlph}{\begin{enumerate}
[label={\normalfont(\makebox[\mywidth]{\alph*})}]}{\end{enumerate}}

\newenvironment{enumerateOurAlphPrim}{\begin{enumerate}
[label={\normalfont(\makebox[\mywidthprim]{\alph*'})}]}{\end{enumerate}}

\newenvironment{enumerateOurAlphPrimPrim}{\begin{enumerate}
[label={\normalfont(\makebox[\mywidthprimprim]{\alph*''})}]}{\end{enumerate}}

\newdimen\mywidth
\sbox0{m}%
\newdimen\mywidthprim
\sbox0{m}%
\newdimen\mywidthprimprim
\sbox0{m}%
\newdimen\mywidthA
\sbox0{m}%
\newdimen\mywidthAprim
\sbox0{m}%
\DeclarePairedDelimiter\set{\{}{\}}

\theoremstyle{plain}
\newtheorem{thm}{Theorem}
\newtheorem*{thm*}{Theorem}

\newtheorem{lemma}[thm]{Lemma}
\newtheorem*{lemma*}{Lemma}

\newtheorem*{cor*}{Corollary}

\newtheorem*{corollary*}{Corollary}

\theoremstyle{remark}

\newtheorem*{problem*}{Open problem}

\newtheorem{claim}[thm]{Claim}

\crefname{obs}{Observation}{Observations}
\theoremstyle{definition}

\newtheorem*{conj*}{Conjecture}
\crefname{lem}{Lemma}{Lemmas}
\crefname{thm}{Theorem}{Theorems}
\crefname{cor}{Corollary}{Corollaries}

\newenvironment{proofclaim}[1][]
	{\par\noindent {\it Proof of Claim~\theclaim}. }{ \strut\hfill$\lozenge$\par\vspace{11pt}}

\newcommand{\tw}{\operatorname{tw}}

\newcommand{\ltw}{\operatorname{ltw}}

\newcommand{\Oh}{\mathcal{O}}

\newcommand{\calA}{\mathcal{A}} 
\newcommand{\calB}{\mathcal{B}} 
\newcommand{\calC}{\mathcal{C}}
\newcommand{\calD}{\mathcal{D}}
 
\newcommand{\calF}{\mathcal{F}}

\newcommand{\calL}{\mathcal{L}}

\newcommand{\calP}{\mathcal{P}} 
\newcommand{\calQ}{\mathcal{Q}} 
\newcommand{\calR}{\mathcal{R}}

\newcommand{\calT}{\mathcal{T}}
\newcommand{\calU}{\mathcal{U}} 
\newcommand{\calV}{\mathcal{V}}
\newcommand{\calW}{\mathcal{W}}

\newcommand{\bigO}{\mathcal{O}}

\newcommand{\lca}{\mathrm{lca}}
\newcommand{\LCA}{\mathrm{LCA}}

\newcommand{\torso}{\mathrm{torso}}

\newcommand{\NN}{\mathbb{N}}

\let\leq\leqslant
\let\geq\geqslant

\let\subset\subseteq

\let\epsilon\varepsilon
\let\phi\varphi

\DeclareMathOperator\WReach{WReach}
\DeclareMathOperator\wcol{wcol}

\DeclareMathOperator\parent{p}
\let\root\relax
\DeclareMathOperator\root{r}

\renewcommand{\setminus}{-}

\renewcommand{\root}{\mathrm{root}}

\title{Centered colorings in minor-closed graph classes}

\begin{document}

\author[Hodor]{Jędrzej Hodor}
\address[J.~Hodor]{Theoretical Computer Science Department, 
Faculty of Mathematics and Computer Science and Doctoral School of Exact and Natural Sciences, Jagiellonian University, Krak\'ow, Poland}
\email{jedrzej.hodor@gmail.com}

\author[La]{Hoang La}
\address[H.~La]{LISN, Universit\'e Paris-Saclay, CNRS, Gif-sur-Yvette, France}
\email{hoang.la.research@gmail.com}

\author[Micek]{Piotr Micek}
\address[P.~Micek]{Theoretical Computer Science Department, 
Faculty of Mathematics and Computer Science, Jagiellonian University, Krak\'ow, Poland}
\email{piotr.micek@uj.edu.pl}

\author[Rambaud]{Clément Rambaud}
\address[C.~Rambaud]{Universit\'e C\^ote d'Azur, CNRS, Inria, I3S, Sophia-Antipolis, France}
\email{clement.rambaud@inria.fr}

\thanks{This research was funded by the National Science Center of Poland under grant UMO-2023/05/Y/ST6/00079 within the WEAVE-UNISONO program. 
Additionally, J.\ Hodor was supported by a Polish Ministry of Education and Science grant (Perły Nauki; PN/01/0265/2022); H.\ La  was supported by ANR PIA funding: ANR-20-IDEES-0002; C.\ Rambaud was supported by ANR DIGRAPHS funding: ANR-19-CE48-0013.}

\begin{abstract}
    A vertex coloring $\phi$ of a graph $G$ is \emph{$p$-centered} if for every connected subgraph $H$ of $G$, either $\phi$ uses more than $p$ colors on $H$, or there is a color that appears exactly once on $H$.
    We prove that for every fixed positive integer $t$, every $K_t$-minor-free graph admits a $p$-centered coloring using 
    $\bigO(p^{t-1})$ colors.
\end{abstract}

\maketitle
\newpage
\section{Introduction}

Let $G$ be a graph, $p$ be a positive integer, and $C$ be a set of colors.
A coloring $\varphi\colon V(G)\to C$ of $G$ is \emph{$p$-centered} if for every connected\footnote{
In this paper, connected graphs are nonnull, that is, they have at least one vertex. A tree is defined as a connected forest, thus, trees and subtrees are also assumed to be nonnull.} subgraph $H$ of $G$, 
either $\varphi$ uses more than $p$ colors on $H$,
or there is a color that appears exactly once on $H$.
The \emph{$p$-centered chromatic number} of $G$, denoted by $\chi_p(G)$, introduced by Nešetřil and Ossona de Mendez~\cite{Nesetril2008}, 
is the least nonnegative integer $k$ such that $G$ admits a $p$-centered coloring using $k$ colors.
The following theorem is the main contribution of this paper.

\begin{thm}\label{thm:main_Kt_minor_free}
    Let $t$ be an integer with $t \geq 2$.
    There exists an integer $c$ such that
    for every $K_t$-minor-free graph $G$ and
    for every positive integer $p$,
    \[
    \chi_p(G) \leq c \cdot p^{t-1}.
    \]
\end{thm}
This improves on the work of~Pilipczuk and Siebertz~\cite{PS19} who proved that for every integer $t$ with $t \geq 2$,
$K_t$-minor-free graphs have $p$-centered chromatic number upper bounded by a polynomial in $p$.
Contrary to \Cref{thm:main_Kt_minor_free}, the degree of their polynomial is not explicitly given and arises from an application of the graph minor structure theorem
by Robertson and Seymour~\cite{GM16}.
On the other hand, Dębski, Micek, Schr\"{o}der, and Felsner~\cite{Dbski2021} showed that there exist $K_t$-minor-free graphs with $p$-centered chromatic number in $\Omega(p^{t-2})$.
Hence, \Cref{thm:main_Kt_minor_free} is tight, up to an $\bigO(p)$ factor. 
Let us provide some context for our theorem.

One of the driving forces in graph theory continues to be the development of efficient algorithms for computationally hard problems for sparse graph classes.
Ne{\v{s}}et{\v{r}}il and Ossona de Mendez~\cite{sparsity} introduced the concepts of bounded expansion and nowhere denseness of classes of graphs.
These notions cover many well-studied classes of graphs, such as 
planar graphs, graphs of bounded treewidth, graphs excluding a fixed minor, graphs of bounded book-thickness, or graphs that admit drawings with a bounded number of crossings per edge.
See~\cite{NOdMW12} by Ne{\v{s}}et{\v{r}}il, Ossona de Mendez, and Wood,
or the recent lecture notes~\cite{notes} of Pilipczuk, Pilipczuk, and Siebertz. 
A key reason why centered colorings have gained 
significant attention is that they encapsulate the concept of bounded expansion.
Indeed, a class of graphs $\calC$ has \emph{bounded expansion} 
if and only if
there is a function $f\colon\mathbb{N}\to\mathbb{N}$ such that 
for all positive integers $p$, we have $\chi_p(G) \leq f(p)$ for every graph $G$ in $\mathcal{C}$, see again~\cite{sparsity}.

Centered colorings are also a crucial tool in designing parameterized algorithms in classes of graphs of 
bounded expansion.
For example,
Pilipczuk and Siebertz~\cite{PS19} showed that if $\calC$ is a class of graphs excluding a fixed minor,
then it can be decided whether a given $p$-vertex graph $H$ is a subgraph of a given $n$-vertex
graph $G$ in $\calC$ in time $2^{\Oh(p\log p)}\cdot n^{\Oh(1)}$ and space $n^{\Oh(1)}$.
This algorithm relies on a general fact that the union of any $p$ color classes in a $p$-centered coloring induces a subgraph of treedepth at most $p$.
Therefore, finding a $p$-centered coloring using $p^{\Oh(1)}$ colors, allows us to reduce the problem to graphs of bounded treedepth, on which the subgraph isomorphism problem can be solved efficiently.
The running times of algorithms based on $p$-centered colorings 
heavily depend on the number of colors used.

Another important family of graph parameters capturing sparsity are
the weak coloring numbers. 
Let $G$ be a graph, let $\sigma$ be a vertex ordering of $G$, and let $r$ be a nonnegative integer.
For all $u$ and $v$ vertices of $G$, we say that
$v$ is \emph{weakly $r$-reachable from $u$ in $(G,\sigma)$}, if there exists
a path between $u$ and $v$ in $G$
containing at most $r$ edges 
such that for every vertex $w$ on the path, $v\leq_{\sigma} w$.
Let $\WReach_r[G, \sigma, u]$ be the set of all vertices that are weakly $r$-reachable from $u$ in $(G,\sigma)$. 
The $r$-\emph{th weak coloring number} of $G$ is defined as
\[\wcol_r(G) = \min_{\sigma}\max_{u \in V(G)}\ |\WReach_r[G, \sigma, u]|,\]
where $\sigma$ ranges over all 
vertex orderings of $G$.

Van den Heuvel, Ossona de Mendez, Quiroz, Rabinovich, and Siebertz~\cite{vdHetal17} established that 
for every integer $t$ with $t\geq2$, for every $K_t$-minor-free graph $G$,
we have $\wcol_r(G)=\Oh(r^{t-1})$.
Their proof hinges on a carefully designed partition of the vertex set of $G$, inspired by Andreae's work~\cite{And86} on the cops and robber game for $K_t$-minor-free graphs.
This partition, now commonly referred to as a chordal partition, mimics an elimination ordering of vertices of a chordal graph. 
Chordal partitions have since become a significant tool in both structural and algorithmic graph theory.
In particular, $\wcol_r(G)=\Oh(r^{t-1})$ is obtained in~\cite{vdHetal17} with a relatively simple and self-contained argument. 
However, there is currently no known analogous partitioning method that works well with $p$-centered colorings. 
In~\cite{PS19}, the bound of $\chi_p(G)=\Oh(p^{g(t)})$ for some function~$g$ was achieved by first proving the result for graphs of bounded Euler genus, and then lifting the argument using the graph minor structure theorem  
by Robertson and Seymour~\cite{GM16}. 
The degree of the polynomial bound in $p$ they obtain depends on the constants of this structure theorem and is not explicitly given.
Consequently, a significant gap exists between the state-of-art upper bounds for weak coloring numbers and the ones for centered chromatic numbers of $K_t$-minor-free graphs. 
In this work, we close this gap. 

The main novelty of the proof is a replacement of 
a Helly-type property well-known for graphs of bounded treewidth (see~\Cref{lemma:helly_property_tree_decomposition}) 
by a statement on centered colorings that works in the broader class of $K_t$-minor-free graphs. See the definition of $(p,c)$-good colorings and \Cref{lemma:Kt_free_graphs_have_centered_Helly_colorings}.
The final centered colorings are constructed in this paper in two steps.
First, we prove this Helly-type property for $K_t$-minor-free graphs following a layered variant of Robertson-Seymour decomposition introduced by Dujmovi\'{c}, Morin, and Wood~\cite{DJMMUW20} (see \Cref{theorem:Kt_free_product_structure_decomposition}).
Second, using this property, we adapt a framework proposed by Illingworth, Scott, and Wood in~\cite{ISW22}.
This in result allows us to move most of the dependency on $t$ from the degree of the final polynomial to a multiplicative constant in front of it.
For a more detailed outline of the proof, see \Cref{sec:outline}.

We remark that we also investigate bounds for $p$-centered chromatic numbers of $X$-minor-free graphs 
for an arbitrary fixed graph $X$.
\cref{thm:main_Kt_minor_free} implies an upper bound in $\Oh(p^{|V(X)|-1})$. 
However, when $X$ is not a clique one might do much better. 
For instance, when $X$ has treedepth $t$, one can show that $\chi_p(G)$ for $X$-minor-free graphs $G$ is in $\Oh(p^{t-1} \log p)$. 
In particular, for any positive integers $s$ and $t$ with $s\leq t$, the maximum $p$-centered chromatic number of $K_{s,t}$-minor-free graphs is in $\Oh(p^s \log p)$. 
A manuscript detailing these results is being prepared.

Graphs with bounded degree have $p$-centered chromatic number in $\Oh(p)$ as proved in~\cite{Dbski2021}.
This, combined with the structure theorem by Grohe and Marx~\cite{Grohe2015}, gives a polynomial in $p$ upper bound on $\chi_p(G)$ for all $G$ excluding $K_t$ as a topological minor. However, the degree of this polynomial is not explicit. We conclude the introduction with the following problem suggestion.
\begin{problem*}
    Let $t$ be an integer with $t \geq 2$.
    Does there exist a nonnegative integer $c$ such that,
    for every $K_t$-topological-minor-free graph $G$ and
    for every positive integer $p$,
    $\chi_p(G) \leq c \cdot p^{t-1}$?
\end{problem*}

\section{Outline}\label{sec:outline}
For a positive integer $k$, we write $[k]=\{1,\ldots,k\}$ and $[0] = \emptyset$.
Let $S$ be a set. 
An \emph{ordering} of a finite set $E$ is a sequence $\sigma = (x_1, \dots, x_{|E|})$ of all the elements of $E$. 
We write $\min_{\sigma} E = x_1$ and $\max_{\sigma} E = x_{|E|}$.
A \emph{coloring} of $S$ is a function $\phi\colon S \rightarrow C$ for some set $C$.
For each $u \in S$, we say that $\varphi(u)$ is the \emph{color} of $u$ 
and subsequently we say that $\varphi(S)$ is the \emph{set of colors used by $\varphi$}.
For two colorings $\phi_1$ and $\phi_2$ of $S$, we define the coloring $\phi_1 \times \phi_2$ of $S$, called the \emph{product coloring} of $\phi_1$ and $\phi_2$, by $(\phi_1\times \phi_2)(u) = (\phi_1(u),\phi_2(u))$ for every $u \in S$.
Given $S' \subset S$, an element $u \in S'$ is a \emph{$\varphi$-center} of $S'$ if the color of $u$ is unique in $S'$ under $\varphi$, in other words, $\varphi(u) \notin \varphi(S' \setminus \{u\})$.
Let $G$ be a graph.
A \emph{coloring} of $G$ is a coloring of $V(G)$.
Recall that a coloring $\varphi$ of $G$ is a $p$-centered coloring of $G$ for a positive integer $p$ if for every connected subgraph $H$ of $G$, either $|\varphi(V(H))| > p$ or $V(H)$ has a $\phi$-center.
A collection $\mathcal{P}$ of subsets of a nonempty set $S$ is a \emph{partition} of $S$ 
if elements of $\mathcal{P}$ are nonempty, pairwise disjoint, and $\bigcup \mathcal{P} = S$.
Given a graph $G$ and a partition $\mathcal{P}$ of $V(G)$, the \emph{quotient graph} $G/\mathcal{P}$ is the graph with the vertex set $\mathcal{P}$ 
and two distinct $P,P' \in \mathcal{P}$ 
are adjacent in $G/\mathcal{P}$ if there are $u \in P$ and $u' \in P'$ such that $uu'$ is an edge in $G$.
The \emph{neighborhood} of $u \in V(G)$ in a graph $G$, denoted by $N_G(u)$, is the set $\{v \in V(G) \mid uv \in E(G)\}$.
For every set $X$ of vertices of a graph $G$, let $N_G(X)=\bigcup_{u \in X} N_G(u) \setminus X$.

For a graph $H$, a \emph{model} of $H$ in a graph $G$ is a family $\big(B_x \mid x \in V(H)\big)$ of disjoint subsets of $V(G)$ such that
\begin{enumerate}
    \item $G[B_x]$ is connected, for every $x \in V(H)$, and
    \item there is an edge between $B_x$ and $B_y$ in $G$, for every $xy \in E(H)$.
\end{enumerate}
If $G$ has a model of $H$, then we say that $H$ is a \emph{minor} of $G$.

A \emph{tree decomposition} of $G$ is a pair $\mathcal{W} = \big(T,(W_x \mid x \in V(T))\big)$
where $T$ is a tree and $W_x \subseteq V(G)$ for every $x \in V(T)$ satisfying the following conditions:
\begin{enumerate}
    \item for every $u \in V(G)$, $T[\{x \in V(T) \mid u \in W_x\}]$ is a connected subgraph of $T$, and
    \item for every edge $uv \in E(G)$, there exists $x \in V(T)$ such that $u,v \in W_x$.
\end{enumerate}
The sets $W_x$ are called the \emph{bags} of $\mathcal{W}$.
The sets $W_x \cap W_y$ for $xy \in E(T)$ are called the \emph{adhesions} of $\mathcal{W}$.
The \emph{width} of $\mathcal{W}$ is $\max_{x \in V(T)} |W_x|-1$,
and the \emph{treewidth} of $G$, denoted by $\tw(G)$, is the minimum width of a tree decomposition of $G$.
An \emph{elimination ordering} of $\calW$\footnote{Equivalently, $(u_1, \dots, u_{|V(G)|})$ is an elimination ordering of $\mathcal{W}$ if and only if it is a perfect elimination ordering of the chordal graph obtained from $G$ by adding all possible edges between vertices in a same bag $W_x$ for every $x \in V(T)$.} is an ordering $(u_1, \dots, u_{|V(G)|})$ of $V(G)$ such that for every $i \in [|V(G)|]$,
there exists $x \in V(T)$ such that 
\[\bigcup \{W_z \mid z \in V(T)\text{ and } u_i \in W_z \} \cap \{u_j \mid j \in [i-1]\} \subseteq W_x.\]

Next, we explore the key ideas underlying the proof of~\Cref{thm:main_Kt_minor_free}.

First, let us discuss the centered colorings of graphs of bounded treewidth.
Let $w$ and $t$ be positive integers with $t \geq 2$, 
and let $G$ be a graph with $\tw(G) \leq w$.
Pilipczuk and Siebertz~\cite{PS19}, showed that $\chi_p(G)\leq \binom{p+w}{w} = 
 \mathcal{O}(p^w)$.
In fact, they proved a stronger statement, namely, they used elimination orderings of tree decompositions to localize potential centers in a structured way.
In order to state this stronger version of the result, let us define a variant of centered colorings in ordered graphs.
Let $\sigma$ be an ordering of $V(G)$.
A coloring $\varphi$ of $G$ is a \emph{$p$-centered coloring} of $(G, \sigma)$ if for every connected subgraph $H$ of $G$, either $|\phi(H)| > p$ or $\min_\sigma V(H)$ is a $\phi$-center of $V(H)$.
We denote by $\chi_p(G, \sigma)$ the least nonnegative integer $k$ such that $(G,\sigma)$ admits a $p$-centered coloring using $k$ colors.
\begin{thm}[\cite{PS19}]\label{thm:Pilipczuk_Siebertz_tw}
    Let $w$ be a positive integer, $G$ be a graph, and $\calW$ be a tree decomposition of $G$ of width at most $w$.
    For every elimination ordering $\sigma$ of $\calW$ and for every positive integer~$p$,
    \[
    \chi_p(G,\sigma) \leq \binom{p+w}{w}.
    \]
\end{thm}

The bound $\chi_p(G) \leq \binom{p+w}{w}$ for graphs with $\tw(G) \leq w$ can be used to obtain a good upper bound on the $p$-centered chromatic number of $K_t$-minor-free graphs of bounded treewidth.
In an influential recent work, 
Illingworth, Scott, and Wood~\cite{ISW22} proved that
every $K_t$-minor-free graph $G$ with $\tw(G) \leq w$
admits a vertex-partition $\calP$ such that $G / \calP$ has treewidth at most $t-2$ and each part of $\calP$ has at most $w+1$ elements\footnote{In other words, there exists a graph $A$ of treewidth at most $t-2$ such that 
$G \subseteq A \boxtimes K_{w+1}$, where $\boxtimes$ denotes the strong product.}. 
This combined with~\Cref{thm:Pilipczuk_Siebertz_tw} implies that for every $K_t$-minor-free graph $G$ of bounded treewidth, we have $\chi_p(G) = \bigO(p^{t-2})$ as we show below.
Let $G$ be a $K_t$-minor-free with $\tw(G) \leq w$ and let 
$\calP$ be a partition of $V(G)$ as above.
Let $\xi$ be a coloring of $G / \calP$ given by~\cref{thm:Pilipczuk_Siebertz_tw} and $\rho$ be a coloring of $G$ using at most $w+1$ colors that is injective on each $P \in \calP$.
For every $P \in \calP$ and $u \in P$, let $\zeta(u) = (\xi(P),\rho(u))$.
To show that $\zeta$ is indeed a $p$-centered coloring of $G$, let $H$ be a connected subgraph of $G$ such that $|\zeta(V(H))| \leq p$.
The parts of $\calP$ that intersect $V(H)$ induce a connected subgraph of $G / \calP$, thus, there is a $\xi$-center $P \in \calP$ of the set of parts of $\calP$ intersecting $V(H)$.
Now, since $\rho$ is injective on $P$, any vertex in $P$ is a $\zeta$-center of $V(H)$.

The previous paragraph sketches the proof of the statement like in~\Cref{thm:main_Kt_minor_free} but for $K_t$-minor-free graphs of bounded treewidth. 
The essential property of the parts in $\calP$ we used is that their sizes are bounded.
There is no hope of obtaining such a partition $\calP$, with each part of size bounded by a constant, in the general case of $K_t$-minor-free graphs.
This forces us to relax the condition on parts, so instead of bounding the maximal size of the part we want to maintain enough structural information.
We accomplish this goal with an extra pre-coloring $\rho$ of vertices of $G$ so that we are able to mimic the proof as in the paragraph above. 
Since $\rho$ uses $\Oh(p)$ colors, 
this ultimately gives the bound $\mathcal{O}(p^{t-1})$ instead of $\mathcal{O}(p^{t-2})$. 
Below is a precise technical statement that we prove later in the paper and a simple argument that the statement together with~\Cref{thm:Pilipczuk_Siebertz_tw} implies~\Cref{thm:main_Kt_minor_free}.

\begin{lemma}\label{lemma:ISW-lifted}
    For every positive integer $t$, there exists a positive integer $c_{\ref{lemma:ISW-lifted}}(t)$ such that for every $K_t$-minor-free graph $G$ and every positive integer $p$, there exists 
    a partition $\mathcal{P}$ of $V(G)$, 
    a tree decomposition $\calW$ of $G/\mathcal{P}$ of width at most $t-2$, 
    an elimination ordering $\sigma = (P_1, \dots, P_\ell)$ of $\calW$, and 
    a coloring $\rho$ of $G$ using at most $c_{\ref{lemma:ISW-lifted}}(t) \cdot (p+1)$ colors 
    such that for every connected subgraph $H$ of $G\left[\bigcup\{Q \in \mathcal{P}\mid Q \geq_\sigma P\}\right]$ with $V(H) \cap P \neq \emptyset$,
    either $|\rho(V(H))| > p$, or $V(H) \cap P$ has a $\rho$-center.
\end{lemma}

\begin{figure}[tp]
  \begin{center}
    \includegraphics{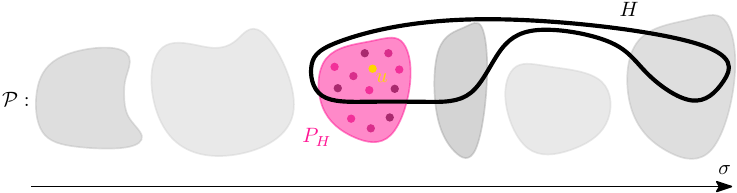}
  \end{center}
  \caption{
    Finding a $\zeta$-center of $V(H)$ in $P_H$. 
  }
  \label{fig:finding-center-in-zeta}
\end{figure}

\begin{proof}[Proof of \Cref{thm:main_Kt_minor_free}]
    Let $c_{\ref{lemma:ISW-lifted}}(t)$ be the constant from~\Cref{lemma:ISW-lifted},
    and let $c = 2^{t-1}\cdot c_{\ref{lemma:ISW-lifted}}(t)$.
    Let $G$ be a $K_t$-minor-free graph,
    and let $p$ be a positive integer.
    Let $\calP$, $\calW$, $\sigma$, and $\rho$ be obtained by~\Cref{lemma:ISW-lifted}.
    Let $\xi$ be a $p$-centered coloring of $G / \calP$ using at most $\binom{p+t-2}{t-2}$ colors obtained from~\Cref{thm:Pilipczuk_Siebertz_tw} applied to $G / \calP$, $\calW$, $\sigma$, and $p$.
    We define a coloring $\zeta$ of $G$ in the following way. 
    For every $P \in \calP$ and $u \in P$,
        \[ \zeta(u) = (\xi(P), \rho(u)).\] 
    Note that 
    \[|\zeta(V(G))| \leq 
    \binom{p+t-2}{t-2} \cdot c_{\ref{lemma:ISW-lifted}}(t) \cdot (p+1) 
    \leq c_{\ref{lemma:ISW-lifted}}(t) \cdot (p+1)^{t-1} \leq  
    c \cdot p^{t-1}.\]
    It remains to show that $\zeta$ is a $p$-centered coloring of $G$.
    Let $H$ be a connected subgraph of $G$ with $|\zeta(V(H))|\leq p$.
    Let $\mathcal{P}_H =  \{P \in \mathcal{P} \mid P \cap V(H) \neq \emptyset\}$ and $P_H = \min_\sigma \mathcal{P}_H$.
    See~\cref{fig:finding-center-in-zeta}.
    Since $|\xi(\mathcal{P}_H)| \leq |\zeta(V(H))| \leq p$ and $\xi$ is a $p$-centered coloring of $(G/\calP,\sigma)$, $P_H$ is a $\xi$-center of $\calP_H$.
    So $\zeta(u) \neq \zeta(v)$ for all $u \in V(H) \cap P_H$ and $v \in V(H) \setminus P_H$.
    Now, $H$ is a connected subgraph of $G\left[\bigcup\{Q \in \mathcal{P}\mid Q \geq_\sigma P_H\}\right]$ with $V(H) \cap P_H \neq \emptyset$ by definition of $P_H$.
    Since $|\rho(V(H))| \leq |\zeta(V(H))| \leq p$, there exists a $\rho$-center $u$ in $V(H) \cap P_H$.
    It follows that $u$ is a $\zeta$-center in $V(H)$.
    We conclude that $\zeta$ is indeed a $p$-centered coloring of $G$.
\end{proof}

The proof of~\Cref{lemma:ISW-lifted} draws on the ideas from~\cite{ISW22}.  
Let $t$ be an integer with $t\geq2$,
and let $G$ be a $K_t$-minor-free graph. 
Next, we describe a key step of the construction (from~\cite{ISW22}) of a partition $\calP$ of vertices of $G$ such that $G/\calP$ has treewidth at most $t-2$ and each part of $\calP$ has at most $\tw(G)+1$ elements. 
Let $R = \{r_1,\dots,r_{t-2}\}$ be a subset of vertices of $G$. 
(One can think of vertices in $R$ as contracted parts of the vertex-partition constructed so far.) 
For simplicity, 
assume that $G-R$ is connected. 
Set $\calF$ to be the family of all connected subgraphs of $G-R$ adjacent to all $r_1,\dots,r_{t-2}$ in $G$. 
Since $G$ is $K_t$-minor-free (and $G-R$ is connected), 
$\calF$ does not contain $t$ pairwise disjoint members, as otherwise, $G$ contains a $K_t$-minor (as shown in~\cref{fig:constructing-Kt-minor}).
Now comes the critical moment when the argument (from~\cite{ISW22}) exploits the fact that $G$ has bounded treewidth and therefore $G$ admits the following Helly-type property.

\begin{lemma}[{\cite[Statement (8.7)]{GM5}}]\label{lemma:helly_property_tree_decomposition}
    For every graph $G$, for every tree decomposition $\mathcal{V}$ of $G$, for every family $\mathcal{F}$ of connected subgraphs of $G$, for every nonnegative integer $d$, either
    \begin{enumerate}
        \item there are $d+1$ pairwise vertex-disjoint subgraphs in $\mathcal{F}$, or
        \item there is a set $Z$ that is the union of at most $d$ bags of $\mathcal{V}$ such that $V(F) \cap Z \neq \emptyset$ for every $F \in \mathcal{F}$. \label{item:helly_property_tree_decomposition:hit}
    \end{enumerate}
\end{lemma}

We take as $\calV$ a tree decomposition of $G$ witnessing treewidth of $G$ and now \Cref{lemma:helly_property_tree_decomposition} implies that there is a set $Z \subset V(G)$ of size at most $t \cdot (\tw(G)+1)$ such that $V(F) \cap Z \neq \emptyset$ for every $F \in \calF$.
This set $Z$ is destined to be a new part in the constructed vertex-partition. 
The key property is that each component $C$ of $G - R - Z$ is not in $\calF$,
and so $C$ is not adjacent to at least one vertex in $R$. 
This allows us to continue the process for $C$ where a non-adjacent to $C$ vertex of $R$ is replaced by the contraction of $Z$.
This sketch omits numerous details, e.g.\ we do not know if $Z$ is connected in $G$.

The key property of a graph $G$ required for the idea in~\cite{ISW22} to work is that for every family $\calF$ of connected subgraphs of $G$, there is a hitting set\footnote{A \emph{hitting set} of a family $\calF$ of connected subgraphs of a graph $G$ is a set $Z \subset V(G)$ such that $V(F) \cap Z \neq \emptyset$ for every $F \in \calF$.} of $\calF$ of size bounded by a function of the packing number\footnote{A \emph{packing number} of a family $\calF$ of connected subgraphs of a graph $G$ is the maximum number of pairwise vertex-disjoint members of $\calF$.} of $\calF$.
Graphs of bounded treewidth indeed satisfy this property as witnessed by~\Cref{lemma:helly_property_tree_decomposition}.
However, $K_t$-minor-free graphs do not admit such a Helly-type property\footnote{Consider the $n \times n$ planar grid $G$ and the family $\calF$ consisting of all the unions of one row and one column in the grid. The packing number of $\calF$ is $1$ but there is no hitting set of $\calF$ that has less than $n$ elements. Additionally, each planar grid is $K_5$-minor-free.}.
The observation that will eventually lead to the proof of \Cref{thm:main_Kt_minor_free} is that we do not need these hitting sets to have bounded size,
as long as we can find a suitable coloring for them.

More precisely, we propose the following technical definition.
Let $p$ and $c$ be positive integers.
A coloring $\phi$ of $G$ is \emph{$(p,c)$-good} if
for every subgraph $G_0$ of $G$,
for every family $\mathcal{F}$ of connected subgraphs of $G_0$,
for every positive integer $d$,
if there are no $d+1$ pairwise disjoint members of $\mathcal{F}$, then
there exists $Z \subseteq V(G_0)$ and $\psi \colon Z \to [c \cdot d]$ such that
\begin{enumerate}[label=(pc\arabic*)]
    \item $V(F) \cap Z \neq \emptyset$ for every $F \in \calF$; \label{centered-Helly-hitting}
    \item for every connected subgraph $H$ of $G_0$ with $V(H) \cap Z \neq \emptyset$, either $|\phi(V(H))|>p$ or $V(H) \cap Z$ has a $(\phi \times \psi)$-center;\label{centered-Helly-coloring} 
    \item for every component $C$ of $G_0-Z$, $N_{G_0}(V(C))$ intersects at most two components of $G_0-V(C)$. \label{centered-Helly-components}
\end{enumerate}

The main technical step of our proof is the following lemma proved in~\cref{sec:building_phi}.

\begin{lemma}\label{lemma:Kt_free_graphs_have_centered_Helly_colorings}
    For every positive integer $t$, there exists a positive integer $c_{\ref{lemma:Kt_free_graphs_have_centered_Helly_colorings}}(t)$ 
    such that
    for every $K_t$-minor-free graph $G$ and every positive integer $p$,
    $G$ admits a $(p,c_{\ref{lemma:Kt_free_graphs_have_centered_Helly_colorings}}(t))$-good coloring using $p+1$ colors.
\end{lemma}

The existence of good colorings in $K_t$-minor-free graphs (\Cref{lemma:Kt_free_graphs_have_centered_Helly_colorings}) suffices to replace the role of~\Cref{lemma:helly_property_tree_decomposition} in the sketched proof framework of~\cite{ISW22} and this way we close the proof of~\Cref{lemma:ISW-lifted}.
See~\Cref{sec:good-to-centered}.

In the remainder of the section, we discuss the ideas behind the proof of~\Cref{lemma:Kt_free_graphs_have_centered_Helly_colorings}.
For simplicity, in this outline, we consider a simpler variant of good colorings, where we drop~\ref{centered-Helly-components}.
Observe that in the definition of good colorings, a subgraph $G_0$ occurs only in~\ref{centered-Helly-components}, hence, if we ignore this item, we can also assume that $G_0 = G$.
Solely for this outline, we say that a coloring is a \emph{simple} $(p,c)$-good coloring if it admits the relaxed definition as discussed above.

For graphs of bounded treewidth, there is a straightforward way of constructing simple good colorings.
Let $G$ be a graph.
We can set $\phi$ to use only one color, and for a family $\calF$ of connected subgraphs of $G$ having no $d+1$ pairwise disjoint members we take $Z$ as given by \Cref{lemma:helly_property_tree_decomposition}, 
so~\ref{centered-Helly-hitting} is satisfied. 
For an injective $\psi\colon Z \to [(\tw(G)+1) \cdot d]$, \ref{centered-Helly-coloring} is clearly satisfied.
This gives a $(p,\tw(G)+1)$-good coloring of $G$ using one color for every positive integer $p$.
A true inspiration for the definition of good colorings comes from the case of graphs of bounded layered treewidth.
Let us first introduce this parameter properly.

A \emph{layered tree decomposition} of $G$ is a pair $(\mathcal{W},\calL)$
where $\mathcal{W}= \big(T,(W_x \mid x \in V(T))\big)$ is a tree decomposition of $G$,
and $\calL=(L_i \mid i \in \NN)$\footnote{By $\NN$ we denote the set of all nonnegative integers.}
is a \emph{layering} of $G$, that is,
a family of pairwise disjoint subsets of $V(G)$ whose union is $V(G)$ and
such that for every edge $uv$ of $G$, there exists $i \in \NN$ such that $u,v \in L_{i} \cup L_{i+1}$.
The \emph{width} of $(\calW,\calL)$ is $\max |W_x \cap L_i|$ over all $x\in V(T)$ and $i \in \NN$.
The \emph{layered treewidth} of $G$, denoted by $\ltw(G)$, is the minimum width of a layered tree decomposition of $G$.

A typical example of a family of graphs of bounded layered treewidth but unbounded treewidth is the family of planar graphs.
Indeed, building upon ideas by Eppstein~\cite{Eppstein1999},
Dujmović, Morin, and Wood~\cite{Dujmovi2017} proved that planar graphs
have layered treewidth at most three.

Next, we describe how to obtain simple good colorings of graphs of bounded layered treewidth.
Let $G$ be a graph, let $p$ be a positive integer, and let $\ltw(G)$ be witnessed by $(\calW,\calL)$, where $\calW=\big(T,(W_x \mid x\in V(T))\big)$ and $\calL=(L_i \mid i \in \NN)$. 
For all $i \in \NN$ and $v\in L_i$, we define 
$\phi(v) = i \bmod (p+1)$. 
We claim that $\phi$ is a simple $(p,\ltw(G))$-good coloring of $G$.
Let $\mathcal{F}$ be a family of connected subgraphs of $G$ with no $d+1$ pairwise disjoint members. 
Let $Z$ be the union of at most $d$ bags of $\calW$ that is a hitting set of $\mathcal{F}$ in $G$ (given by~\Cref{lemma:helly_property_tree_decomposition}). 
In particular, $|Z \cap L_i| \leq \ltw(G) \cdot d$ for every $i\in \NN$.
Let $\psi \colon Z \to [\ltw(G) \cdot d]$ be any coloring that is injective on $Z \cap L_i$ for every $i \in \NN$.
It remains to verify~\ref{centered-Helly-coloring}.
Let $H$ be a connected subgraph of $G$ with $V(H)\cap Z\neq\emptyset$. 
Since $H$ is connected and $\calL$ is a layering of $G$, $V(H)$ intersects a set of consecutive layers in $\calL$.
If this set has at least $p+1$ elements, then $|\varphi(V(H))| > p$, and otherwise any element of $V(H) \cap Z$ is a $\phi\times\psi$-center of $V(H) \cap Z$.

In order to lift this idea to the class of general $K_t$-minor-free graphs, we use a graph decomposition introduced by Dujmović, Morin, and Wood~\cite{Dujmovi2017}.
We first formally define this decomposition and then discuss how it can be used to construct good colorings.

For a graph $G$ and a subset $W$ of vertices of $G$, the \emph{torso} of $W$ in $G$, denoted by $\torso_G(W)$, is the graph with the vertex set $W$ and two vertices are adjacent if they can be connected in $G$ by a path whose internal vertices do not belong to $W$.
For a graph $G$ and a positive integer $c$, a \emph{layered Robertson-Seymour decomposition} of $G$ (\emph{layered RS-decomposition} for short) of width at most $c$ 
is a tuple 
    \[(T,\calW,\calA,\calD,\calL)\]
where $T$ is a tree and
\begin{enumerate}[label=(lrs\arabic*)]
    \item $\calW = \big(T,(W_x \mid x \in V(T))\big)$ is a tree decomposition of $G$ where adhesions have size at most $c$; \label{LRS:adhesion}
    \item $\calA = \big(A_x \mid x \in V(T)\big)$ where $A_x \subset W_x$ and $|A_x| \leq c$ for every $x \in V(T)$; \label{LRS:apices}
    \item $\calD = \big(\calD_x \mid x \in V(T)\big)$ where $\calD_x = \big(T_x,(D_{x,z} \mid z \in V(T_x))\big)$ is a tree decomposition of $\torso_G(W_x)-A_x$ for every $x \in V(T)$; \label{LRS:tree_decompositions}
    \item $\calL = \big(\calL_x \mid x \in V(T)\big)$ where $\calL_x = (L_{x,i}\mid i \in \mathbb{N})$ is a layering of $\torso_G(W_x)-A_x$ for every $x\in V(T)$; \label{LRS:layering}
    \item $|D_{x,z} \cap L_{x,i}| \leq c$ for all $x \in V(T)$, $z \in V(T_x)$, and $i \in \NN$.  \label{LRS:ltw}
\end{enumerate}

\begin{figure}[tp]
  \begin{center}
    \includegraphics{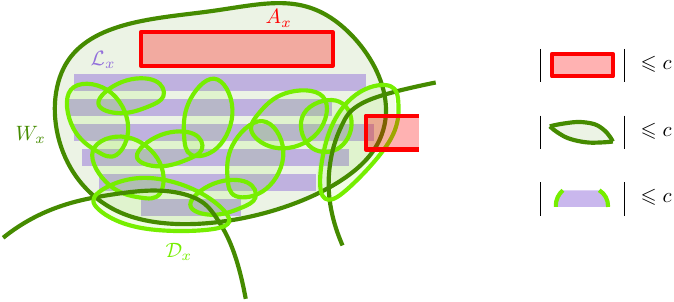}
  \end{center}
  \caption{
    A bag $W_x$ of $\calW$ in a layered RS-decomposition $(\calT, \calW, \calA, \calD, \calL)$ of width at most $c$. 
    The set $A_x$ (in red) is included in $W_x$. 
    The graph $\torso_G(W_x)-A_x$ has a layering $\calL_x$ (in purple) and a tree decomposition $\calD_x$ (in green).
    Note that for every $y \in V(T) - \{x\}$, $W_x \cap W_y - A_x$ is a clique in $\torso_G(W_x) - A_x$, and so, it is contained in a single bag of $\calD_x$ and in at most two layers of $\calL_x$. 
  }
  \label{fig:lrs-outline}
\end{figure}

See~\Cref{fig:lrs-outline}.

\begin{thm}[{\cite[Theorem~22 and Lemma~26]{Dujmovi2017}}]\label{theorem:Kt_free_product_structure_decomposition}
    For every positive integer $t$, there is a positive integer $c_{\ref{theorem:Kt_free_product_structure_decomposition}}(t)$ 
    such that every $K_t$-minor-free graph $G$ admits a layered RS-decomposition of width at most~$c_{\ref{theorem:Kt_free_product_structure_decomposition}}(t)$.
\end{thm}

Given a layered RS-decomposition of a graph $G$ of bounded width, the strategy to find a simple good coloring of $G$ is the following.
Let $p$ be a positive integer.
The intention is to mimic the proof for the bounded layered treewidth case.
That is, ignoring the overlapping between the bags of $\calW$, we set $\varphi(v) \equiv i \bmod (p+1)$ for every $i \in \NN$, and for every $v \in L_{x,i}$.
We claim that $\varphi$ is a simple $(p,c)$-good coloring, where $c$ depends only on the width of the given decomposition.
Let $\mathcal{F}$ be a family of connected subgraphs of $G$ with no $d+1$ pairwise disjoint members.
We glue all the tree decompositions in $\calD$ in a natural way obtaining a tree decomposition $\calU$ of $G$.
Note that the vertices of $A_x$ do not occur in bags of tree decompositions of $\calD$, but we can just add them to all bags of $\calU$ corresponding to $W_x$.
The important property of $\calU$ is that~\ref{LRS:ltw} is preserved.
Next, we apply~\Cref{lemma:helly_property_tree_decomposition} to $\calF$ and $\calU$ obtaining a hitting set $Z$ of $\mathcal{F}$ in $G$, which is the union of at most $d$ bags of $\calU$.
Note that there are also at most $d$ bags of $\calW$, whose union contains $Z$.
Say that these bags correspond to the vertices in $X \subset V(T)$.
The next step is to ``disconnect'' elements of $Z$ in different bags of $\calW$.
To this end, we root $T$ arbitrarily and we define $B$ as the union of $\bigcup_{x \in X}A_x$ and the union of adhesions between $W_x$ and the bag corresponding to the parent of $x$ in $\calW$ for all $x \in X$.
In particular, the size of $B$ is bounded.
Finally, we define a coloring $\psi$ of $Z$ so that $\psi$ is injective on $B$ and $\psi$ is injective on $L_{x,i} \cap Z$ but the colors used are disjoint from $\psi(B)$.
Checking that $\psi$ witnesses $\varphi$ being a simple $(p,c)$-good coloring of $G$ is very similar to the bounded layered treewidth case.

\section{Preliminaries} \label{sec:prelimiaries}

Let $T$ be a tree rooted in a vertex $r$. 
For every vertex $v$ of $T$ with $v\neq r$, 
we define $\parent(T,v)$ to be the vertex following $v$ on the path from $v$ to $r$ in $T$.
For a subtree $T'$ of $T$, we define $\root(T')$ as the unique vertex in $T'$ closest to $r$ in $T$.

Let $u$, $v$ be two (not necessarily distinct) vertices of $T$.
The \emph{lowest common ancestor} of $u$ and $v$ in $T$, denoted by $\lca(T,u,v)$, is the furthest vertex from the root that has $u$ and $v$ as descendants.
Let $Y \subset V(T)$.
The \emph{lowest common ancestor closure} of $Y$ in $T$ is the set $\LCA(T,Y)=\{\lca(T,u,v)\mid u,v\in Y\}$.  
Observe that $\LCA(T,\LCA(T,Y)) = \LCA(T,Y)$.
The following two lemmas are folklore. See e.g.~\cite[Lemma~8]{DHHJLMMRW24} for a proof.

\begin{lemma}\label{lemma:increase_X_in_a_tree}
    Let $m$ be a positive integer.
    Let $T$ be a tree,
    and let $Y$ be a set of $m$ vertices of $T$.
    Let $X=\LCA(T,Y)$.
    Then, $|X|\leq 2m-1$ and for every component $C$ of $T-X$, $|N_T(V(C))| \leq 2$.
\end{lemma}

For a tree $T$ and $xy \in E(T)$, we denote by $T_{x|y}$ the component of $x$ in $T \setminus \{xy\}$.
Note that $T \setminus \{xy\}$ is the disjoint union of $T_{x|y}$ and $T_{y|x}$.
A tree decomposition
$\big(T,(W_x\mid x\in V(T))\big)$ of a graph $G$ is \emph{natural} in $G$ if 
\[\textstyle\text{for every edge $xy$ in $T$, 
the graph $G\left[\bigcup_{z\in V(T_{x|y})} W_z\right] \textrm{is connected.}$}\]

\begin{lemma}\label{lemma:increase_X_to_have_small_interfaces}
    Let $m$ be a positive integer.
    Let $G$ be a graph,
    and let $\mathcal{W}=\big(T,(W_x\mid x\in V(T))\big)$ 
    be a tree decomposition of $G$.
    Let $X$ be a set of $m$ vertices of $T$.
    Then $Z=\bigcup_{x\in\LCA(T,X)}W_x$ is the union of at most $2m-1$ bags of $\mathcal{W}$ 
    such that for every component $C$ of $G-Z$, $N_G(V(C))$
    is a subset of the union of at most two bags of $\mathcal{W}$.
    Moreover, if $\mathcal{W}$ is natural, then $N_G(V(C))$ intersects at most two components of $G-V(C)$.
\end{lemma}

We will also need the following simple combinatorial fact.
Let $T$ be a tree and let $\calQ$ a collection of connected subgraphs of $T$ whose vertex sets partition $V(T)$.
For every $X \subset V(T)$, let $\calQ(X) = \{Q \in \calQ \mid V(Q) \cap X \neq \emptyset\}$.

   \begin{lemma}\label{lemma:LCA-stuff}
        Let $T$ be a tree and let $\calQ$ be a collection of connected subgraphs of $T$ whose vertex sets partition $V(T)$.
        Let $X, Y \subseteq V(T)$ with $X\subseteq Y$ and $\LCA(T,X) = X$.
        If $\calQ(X) = \calQ(Y)$,
        then $\calQ(Y)=\calQ(\LCA(T,Y))$.
    \end{lemma}

    \begin{proof}
        Suppose that $\calQ(X)=\calQ(Y)$.
        Since $Y \subset \LCA(T,Y)$, we have $\calQ(Y) \subset \calQ(\LCA(T,Y))$.
        Thus, it suffices to prove that $\calQ(\LCA(T,Y)) \subset \calQ(Y)$.
        Consider $Q \in \calQ(\LCA(T,Y))$.
        There exist $y_1,y_2 \in Y$ such that $\lca(T,y_1,y_2) \in Q$.
        Let $Q_1,Q_2 \in \mathcal{Q}$ be such that
        $y_1 \in V(Q_1)$, $y_2 \in V(Q_2)$.
        In particular, $Q_1,Q_2 \in \calQ(Y) = \calQ(X)$.
        Hence, there exist $x_1 \in V(Q_1) \cap X$ and $x_2 \in V(Q_2) \cap X$.
        
        If the roots of $Q_1$ and $Q_2$ are not in an ancestor-descendant relation in $T$, then for all $q_1\in V(Q_1)$ and $q_2 \in V(Q_2)$, we have $\lca(T,q_1,q_2) = \lca(T, \root(Q_1),\root(Q_2))$.
        In particular, 
            \[\lca(T,x_1,x_2) = \lca(T, \root(Q_1),\root(Q_2)) = \lca(T,y_1,y_2) \in V(Q).\]
        Observe that $\lca(T,x_1,x_2) \in \LCA(T,X) = X$, thus, $V(Q) \cap X \neq \emptyset$, and so, $Q \in \calQ(X) = \calQ(Y)$. 
        
        Therefore, we can assume that one of the roots say $\root(Q_1)$, is an ancestor of the other $\root(Q_2)$ in $S$.
        It follows that for all $q_1\in V(Q_1)$ and $q_2 \in V(Q_2)$, $\lca(T,q_1,q_2)$ is either equal to $\root(Q_1)$ or is a descendant of $\root(Q_1)$.
        In both cases, $\lca(T,q_1,q_2)$ lies in the path from $q_1$ to $\root(Q_1)$ in $S$.
        Since $Q_1$ is connected, we obtain that $\lca(T,q_1,q_2) \in V(Q_1)$.
        In particular, $\lca(T,x_1,x_2) \in V(Q_1)$, and so, $Q_1 = Q$ implying $Q \in \calQ(Y)$, which ends the proof.
    \end{proof}
    
Finally, we need a simple fact about torsos.
We give the proof for completeness.

\begin{lemma}\label{lemma:projections_on_a_torso_stay_connected}
    Let $G$ be a graph, and let $W \subseteq V(G)$.
    If $H$ is a connected subgraph of $G$ intersecting $W$, then
    $V(H) \cap W$ induces a connected subgraph of $\torso_G(W)$.
\end{lemma}
\begin{proof}
    Let $H$ be a connected subgraph of $G$ intersecting $W$, and suppose to the contrary that the subgraph $H_W$ of $\torso_G(W)$ induced by $V(H) \cap W$ is not connected.
    Let $u,v \in V(H_W)$ be vertices in distinct components of $H_W$ such that the distance between $u$ and $v$ is minimal in $H$.
    Note that the internal vertices of a shortest path between $u$ and $v$ in $H$ do not lie in $W$ as otherwise we obtain a pair of vertices of $H_W$ in distinct components that are closer in $H$.
    It follows that $u$ and $v$ are adjacent in $\torso_G(W)$, and so, in $H_W$, which contradicts the assumption and completes the proof.
\end{proof}

\section{From good colorings to centered colorings}\label{sec:good-to-centered}

In this section, we prove~\Cref{lemma:ISW-lifted} assuming~\Cref{lemma:Kt_free_graphs_have_centered_Helly_colorings} (proved in~\cref{sec:building_phi}).
The lemma below is an inductive setup of the proof and is inspired by the framework given in~\cite{ISW22}.

\begin{lemma}\label{thm:technical_Kt_minor_free_partition}
    Let $t$ be an integer with $t \geq 2$, and
    let $p,c$ be positive integers.
    Let $G$ be a $K_t$-minor-free graph, let $r$ be an integer with $0 \leq r \leq t-2$, and let $R_1, \dots, R_r$ be pairwise disjoint subsets of $V(G)$ with $|R_i| \in \{1,2\}$ for every $i \in [r]$.
    Let $\phi$ be a $(p,c)$-good coloring of $G - \bigcup_{i \in [r]} R_i$.
    There are a partition $\mathcal{P}$ of $V(G)$, a tree decomposition $\calW = \big(T,(W_x \mid x \in V(T))\big)$ of $G/\mathcal{P}$ of width at most $t-2$,
    and an elimination ordering $\sigma = (P_1, \dots, P_\ell)$ of $\calW$ such that
    \begin{enumerateOurAlph}
        \item $P_i = R_i$ for every $i \in [r]$; \label{item:thm:technical_Kt_minor_free_partition:Pi=Ri}
        \item there exists $s \in V(T)$ such that $R_1, \dots, R_r \in W_s$; \label{item:thm:technical_Kt_minor_free_partition:Ri_in_Ws}
        \item for every $P \in \mathcal{P} \setminus \{R_1, \dots, R_r\}$, \label{item:thm:technical_Kt_minor_free_partition:parts_have_small_pcentered_colorings}
            there is a coloring $\psi_P\colon P \to [c \cdot 2^{t-2}(t-1)]$ such that
            for every connected subgraph $H$ of $G\left[\bigcup\{Q \in \mathcal{P}\mid Q \geq_\sigma P\}\right]$ with $V(H) \cap P \neq \emptyset$,
            either $|\phi(V(H))| > p$, or $V(H) \cap P$ has a $\phi \times \psi_P$-center.
    \end{enumerateOurAlph}
\end{lemma}

\begin{proof}
    Let $\mathcal{R} = \{R_1, \dots, R_r\}$ and $R = \bigcup_{i\in [r]} R_i$.
    We proceed by induction on $(|V(G-R)|,r)$ in the lexicographic order. 
    For the base case, assume that $V(G-R) = \emptyset$. 
    Then we set $\mathcal{P}=\mathcal{R}$ and $\sigma = (R_1, \dots, R_r)$.
    Let $T=K_1$ and $W_s = \mathcal{P}$ for $s$ the unique vertex of $T$.
    Then \ref{item:thm:technical_Kt_minor_free_partition:Pi=Ri} and~\ref{item:thm:technical_Kt_minor_free_partition:Ri_in_Ws} hold by construction,
    and~\ref{item:thm:technical_Kt_minor_free_partition:parts_have_small_pcentered_colorings} is vacuously true.
    From now on, we assume $V(G-R) \neq \emptyset$.

    If $r < t-2$, then we set $R_{r+1}$ to be an arbitrary singleton in $V(G-R)$,
    and we apply induction on $G,(R_1, \dots, R_{r+1}),\phi\vert_{V(G-R-R_{r+1})}$.
    This gives a partition $\calP$ of $V(G)$, a tree decomposition $\calW$ of $G/\mathcal{P}$,
    and an elimination ordering $\sigma$ of $\calW$ satisfying all three items.
    Items~\ref{item:thm:technical_Kt_minor_free_partition:Pi=Ri} and~\ref{item:thm:technical_Kt_minor_free_partition:Ri_in_Ws} are now clear by induction.
    Item~\ref{item:thm:technical_Kt_minor_free_partition:parts_have_small_pcentered_colorings} is clear for all $P \in \mathcal{P} \setminus \{R_1, \dots, R_r,R_{r+1}\}$, thus, it suffices to argue it for $P = R_{r+1}$.
    Recall that $R_{r+1} = \{v\}$ for some $v \in V(G)$.
    We set $\phi_{R_{r+1}}(v) = 1$.
    Observe that in this case for every subgraph $H$ of $G$ intersecting $R_{r+1}$, $V(H) \cap R_{r+1}$ has a $\varphi \times \psi_{R_{r+1}}$-center. 
    This completes the proof in the case $r<t-2$. 
    Therefore, from now on, we assume $r=t-2$.

    Next, suppose that $G - R$ is not connected, and let $\mathcal{C}$ be the family of all the components of $G-R$.
    Fix some $C \in \mathcal{C}$, and
    let $G_C = G[V(C) \cup R]$.
    By induction hypothesis applied to $G_C$, $(R_1, \dots, R_{t-2})$, and $\phi\vert_{V(C)}$,
    there are a partition $\mathcal{P}_C$ of $V(G_C)$, a tree decomposition $\mathcal{W}_C = \big(T_C,(W_{C,x} \mid x \in V(T_C))\big)$ of $G_C/\mathcal{P}_C$
    of width at most $t-2$
    and an elimination ordering $\sigma_C = (P_{C,1}, \dots, P_{C,\ell_C})$ of $\mathcal{W}_C$ such that
    \begin{enumerateOurAlphPrim}
        \item $P_{C,i} = R_i$ 
        for every $i \in [r]$; 
        \item there exists $s_C \in V(T_C)$ such that $R_1, \dots, R_{t-2} \in W_{C,s_C}$; 
        \item for every $P \in \mathcal{P}_C \setminus \{R_1, \dots, R_{t-2}\}$, \label{item:thm:technical_Kt_minor_free_partition:parts_have_small_pcentered_colorings'}
            there is a coloring $\psi_{C,P}\colon P \to [c \cdot 2^{t-2}(t-1)]$ such that
            for every connected subgraph $H$ of $G_C\left[\bigcup \{Q \in \mathcal{P}_C \mid Q \geq_{\sigma_C} P\}\right]$ with $V(H) \cap P \neq \emptyset$,
            either $|\phi(V(H))| > p$, or $V(H) \cap P$ has a $\phi \times \psi_{C,P}$-center.
    \end{enumerateOurAlphPrim}

    Then, let $\mathcal{P} = \bigcup_{C \in \mathcal{C}} \mathcal{P}_C$,
    and
    let $\sigma$ be the concatenation of 
    $(R_1, \dots, R_{t-2})$, $(P_{C_1,t-1}, \dots, P_{C_1,\ell_C}), \dots,(P_{C_a,t-1}, \dots, P_{C_a,\ell_C})$
    for an arbitrary ordering $(C_1, \dots, C_a)$ of $\mathcal{C}$.
    Let $T$ be obtained from the disjoint union of $T_C$ over all $C \in \mathcal{C}$ by adding a new vertex $s$ with the
    neighborhood $\{s_C \mid C \in \mathcal{C}\}$.
    For every $C \in \mathcal{C}$ and every $x \in V(T_C)$, let $W_x = W_{C,x}$.
    Additionally, let $W_s = \{R_1, \dots, R_{t-2}\}$.
    By construction, $\mathcal{W} = \big(T,(W_x \mid x \in V(T))\big)$ is a tree decomposition of $G/\mathcal{P}$ of width at most $t-2$
    and $\sigma$ is an elimination ordering of $\mathcal{W}$.
    Moreover, \ref{item:thm:technical_Kt_minor_free_partition:Pi=Ri} and~\ref{item:thm:technical_Kt_minor_free_partition:Ri_in_Ws} are clearly satisfied.
    It remains to show that~\ref{item:thm:technical_Kt_minor_free_partition:parts_have_small_pcentered_colorings} holds.
    Let $P \in \mathcal{P} \setminus \{R_1, \dots, R_{t-2}\}$.
    There exists $C \in \mathcal{C}$ such that $P \in \mathcal{P}_C$.
    Let $\psi_P = \psi_{C,P}$.
    Let $H$ be a connected subgraph of $G\left[\bigcup\{Q \in \mathcal{P}\mid Q \geq_\sigma P\}\right]$ with $V(H) \cap P \neq \emptyset$.
    Note that $H$ is a connected subgraph of $G_C\left[\bigcup\{Q \in \mathcal{P}_C\mid Q \geq_{\sigma_C} P\}\right]$.
    Hence by~\ref{item:thm:technical_Kt_minor_free_partition:parts_have_small_pcentered_colorings'} either $|\phi(V(H))|>p$ or $V(H)\cap P$ has a $(\phi\times\psi_{C,P})$-center $u$.
    In the former case, we are immediately satisfied, and in the latter case, $u$ is also a $(\phi\times\psi_P)$-center of $V(H)\cap P$.
    Therefore,~\ref{item:thm:technical_Kt_minor_free_partition:parts_have_small_pcentered_colorings}  holds, which concludes the case where $G-R$ is not connected.
    From now on, we assume that $G-R$ is connected.

    Suppose that there exists $i \in [t-2]$ such that $N_G(R_i) \cap V(G-R) = \emptyset$.
    In this case, we call induction on $G-R_i, (R_1, \dots, R_{i-1},R_{i+1},\dots,R_{t-2}), \phi$.
    We deduce that
    there is a partition $\mathcal{P'}$ of $V(G-R_i)$, a tree decomposition $\calW' = \big(T',(W_x \mid x \in V(T))\big)$ of $(G-R_i)/\mathcal{P'}$ of width at most $t-2$,
    and an elimination ordering $\sigma' = (P'_1, \dots, P'_{\ell'})$ of $\calW'$ satisfying all three items.
    In particular there exists $s' \in V(T')$ such that $R_1, \dots, R_{i-1},R_{i+1},\dots, R_{t-2} \in W_{s'}$.
    Let $\mathcal{P}=\calP'\cup\set{R_i}$,
    and let $\sigma = (P'_1, \dots, P'_{i-1},R_i,P'_{i+1}, \dots, P'_{\ell'})$.
    Let $T$ be obtained from $T'$ by adding a new neighbor $s$ of $s'$.
    Finally, let $W_s = \{R_1, \dots, R_{t-2}\}$ and $W_x = W_x'$ for every $x \in V(T')$.
    We obtain that $\mathcal{W} = \big(T,(W_x \mid x \in V(T))\big)$ is a tree decomposition of $G/\mathcal{P}$
    of width at most $t-2$ and $\sigma$ is an elimination ordering of $\mathcal{W}$.
    It is immediate from the induction hypothesis 
    that $\mathcal{P}$, $\calW$, and $\sigma$ satisfy~\ref{item:thm:technical_Kt_minor_free_partition:Pi=Ri}-\ref{item:thm:technical_Kt_minor_free_partition:parts_have_small_pcentered_colorings}. 
    From now on, we assume that for every $i \in [t-2]$, there is an edge between $V(G-R)$ and $R_i$.

    This completes the series of simple reductions.
    As a result we can assume the following: $V(G - R) \neq \emptyset$, $r = t-2$, $G - R$ is connected, and for every $i \in [t-2]$, there is an edge between $V(G-R)$ and $R_i$. 
    We proceed with the main argument. 

\begin{figure}[tp]
  \begin{center}
    \includegraphics{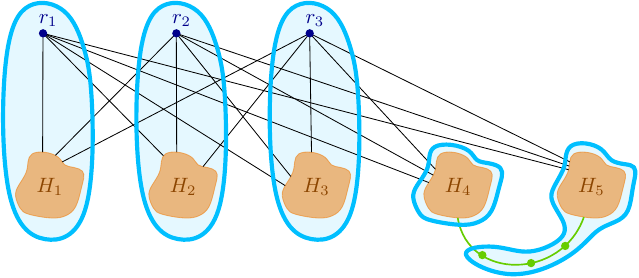}
  \end{center}
  \caption{
      After pigeonholing pairwise disjoint members of $\calF$, we obtain a situation as in the figure.
      Here, $t = 5$.
      The model of $K_t$ is constructed in blue.
  }
  \label{fig:constructing-Kt-minor}
\end{figure}

    Let $\mathcal{F}$ be the family of all connected subgraphs $H$ of $G-R$ such that $N_G(R_i) \cap V(H) \neq \emptyset$ for every $i \in [t-2]$.
    We claim that there are no $2^{t-2}(t-1) + 1$ pairwise disjoint members of $\mathcal{F}$.
    Suppose to the contrary that there are $2^{t-2}(t-1) + 1$ pairwise disjoint members of $\mathcal{F}$.
    Since for each $i \in [t-2]$, $|R_i| \leq 2$, and every member of $\mathcal{F}$ has a neighbor in $R_i$,
    by the pigeonhole principle, there exist pairwise disjoint $H_1,\dots,H_t \in \mathcal{F}$
    and there exists $r_i\in R_i$ for each $i\in[t-2]$ such that for all $j\in[t]$ and $i\in[t-2]$,
    we have that $H_j$ contains a neighbor of $r_i$ in $G$.
    Since $G-R$ is connected, there exist distinct $i_1,i_2 \in [t]$ such that there is a path between $V(H_{i_1})$ and $V(H_{i_2})$ in $G-R$
    internally disjoint from $\bigcup_{j \in [t]} V(H_j)$.
    Let $A$ be the set of internal vertices of this path.
    Without loss of generality, suppose that $i_1 = t-1$ and $i_2=t$.
    Then 
        \[\big(V(H_1) \cup \{r_1\}, \dots, V(H_{t-2}) \cup \{r_{t-2}\}, V(H_{t-1}), V(H_t)\cup A\big)\] 
    is a model of $K_t$ in $G$ (see~\cref{fig:constructing-Kt-minor}), which contradicts the assumption as $G$ is $K_t$-minor-free.
    This contradiction concludes the proof that $\calF$ has no $2^{t-2}(t-1)+1$ pairwise disjoint members.

    Since $\phi$ is a $(p,c)$-good coloring of $G-R$, it follows that
    there exists a set $Z \subseteq V(G-R)$ and $\psi_Z \colon Z \to [c \cdot 2^{t-2}(t-1)]$ such that
    \begin{enumerate}[label=(pc\arabic*')]
        \item $V(F) \cap Z \neq \emptyset$ for every $F \in \mathcal{F}$; \label{items:pc-coloring-local-thm1-hitting}
        \item for every connected subgraph $H$ of $G-R$ with $V(H) \cap Z \neq \emptyset$, either $|\phi(V(H))|>p$, or $V(H) \cap Z$ has a $\phi \times \psi_Z$-center; \label{items:pc-coloring-local-thm1-centers}
        \item for every component $C$ of $G-R-Z$, $N_{G-R}(V(C))$ intersects at most two components of $G-R-V(C)$. \label{items:pc-coloring-local-thm1-components}
    \end{enumerate}
    Note that $G-R \in \mathcal{F}$, and so $\mathcal{F}$ is nonempty.
    In particular, $Z$ is nonempty.
    Let $\mathcal{C}$ be the family of all the components of $G-R-Z$.

\begin{figure}[tp]
  \begin{center}
    \includegraphics{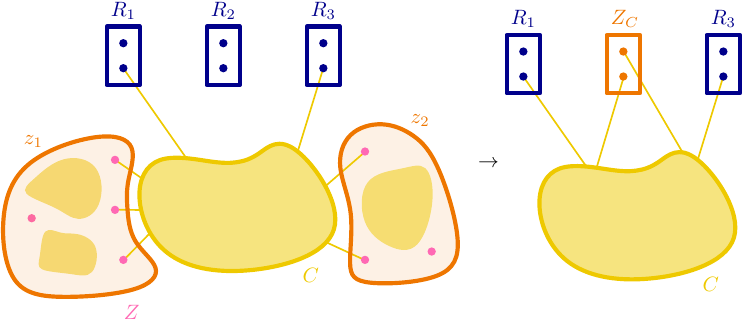}
  \end{center}
  \caption{
    The pink vertices depict the set $Z$.
    The yellow pieces are components of $G - R - Z$. 
    In the figure, $i=2$, that is, $C$ does not have a neighbor in $R_2$.
    The set $Z_C$ consists of two vertices $z_1$ and $z_2$ obtained by contracting the orange parts.
    In order to apply induction, we replace $R_2$ with $\{z_1,z_2\}$. 
  }
  \label{fig:contracting-the-rest}
\end{figure}
    
    Fix some $C \in \mathcal{C}$.
    Since $V(C) \cap Z = \emptyset$, by~\ref{items:pc-coloring-local-thm1-hitting}, we have, $C \notin \mathcal{F}$.
    In particular, there exists $i \in [t-2]$ such that $N_G(R_i) \cap V(C) = \emptyset$.
    Let $G_C$ be obtained from $G - R_i$ by contracting each component of $G-R-V(C)$ into a single vertex.
    See~\cref{fig:contracting-the-rest}.
    Let $Z_C$ be the set of the resulting contracted vertices.
    Since $Z \neq \emptyset$, we have $1 \leq |Z_C|$, and since $G-R$ is connected, by~\ref{items:pc-coloring-local-thm1-components}, $|Z_C| \leq 2$.
    Moreover, $G_C$ is a minor of $G$, thus, $G_C$ is $K_t$-minor-free.
    Let $\mathcal{R}_C = \mathcal{R} \setminus \{R_i\} \cup \{Z_C\}$, and let $(R_{C,1}, \dots, R_{C,t-2})$ be an arbitrary ordering of $\mathcal{R}_C$.
    Since $Z \neq \emptyset$, we have 
    \[|V(G_C-\textstyle\bigcup \mathcal{R}_C)|<|V(G-R)|.\]
    Hence by induction hypothesis applied to $G_C$, $(R_{C,1}, \dots, R_{C,t-2})$, and $\phi\vert_{V(C)}$, 
    there is a partition $\mathcal{P}_C$ of $V(G_C)$, a tree decomposition $\mathcal{W}_C = \big(T_C,(W_{C,x} \mid x \in V(T_C))\big)$ of $G_C/\mathcal{P}_C$
    of width at most $t-2$
    and an elimination ordering $\sigma_C = (P_{C,1}, \dots, P_{C,\ell_C})$ of $\mathcal{W}_C$ such that
    \begin{enumerateOurAlphPrimPrim}
        \item $P_{C,i} = R_{C,i}$ for every $i \in [t-2]$; 
        \item there exists $s_C \in V(T_C)$ such that $R_{C,1}, \dots, R_{C,r} \in W_{C,s_C}$; and 
        \item for every $P \in \mathcal{P}_C \setminus \{R_{C,1}, \dots, R_{C,t-2}\}$, 
        there is a coloring $\psi_{C,P}\colon P \to [c \cdot 2^{t-2}(t-1)]$ such that
        for every connected subgraph $H$ of $G_C\left[\bigcup\{Q \in \mathcal{P}_C\mid Q \geq_{\sigma_C} P\}\right]$ with $V(H) \cap P \neq \emptyset$,
        either $|\phi(V(H))| > p$, or $V(H) \cap P$ has a $\phi \times \psi_{C,P}$-center.
        \label{item:thm:technical_Kt_minor_free_partition:parts_have_small_pcentered_colorings''}
    \end{enumerateOurAlphPrimPrim}

    Then let $\mathcal{P} = \bigcup_{C \in \mathcal{C}} (\mathcal{P}_C \setminus \mathcal{R}_C) \cup \mathcal{R} \cup \{Z\}$,
    and let $\sigma$ be the concatenation of $(R_1, \dots, R_{t-2}, Z),
    (P_{C_1,t-1}, \dots, P_{C_1,\ell_{C_1}}), \dots, (P_{C_a,t-1}, \dots, P_{C_a,\ell_{C_a}})$
    for an arbitrary ordering $(C_1, \dots, C_a)$ of $\mathcal{C}$.
    Let $T$ be obtained from the disjoint union of $T_C$ over all $C \in \mathcal{C}$ by adding a new vertex $s$ with the
    neighborhood $\{s_C \mid C \in \mathcal{C}\}$.
    For every $C \in \mathcal{C}$ and every $x \in V(T_C)$, let $W_x = W_{C,x}$.
    Additionally, let $W_s = \{R_1, \dots, R_{t-2},Z\}$.
    By construction, $\mathcal{W} = \big(T,(W_x \mid x \in V(T))\big)$ is a tree decomposition of $G/\mathcal{P}$ of width at most $t-2$
    and $\sigma$ is an elimination ordering of $\mathcal{W}$.
    Moreover, \ref{item:thm:technical_Kt_minor_free_partition:Pi=Ri} and~\ref{item:thm:technical_Kt_minor_free_partition:Ri_in_Ws} are clearly satisfied.
    
    It remains to show that~\ref{item:thm:technical_Kt_minor_free_partition:parts_have_small_pcentered_colorings} holds.
    Let $P \in \mathcal{P} \setminus \{R_1,\dots, R_{t-2}\}$.
    First, suppose that $P \neq Z$. 
    Then there exists $C \in \mathcal{C}$ such that $P \in \mathcal{P}_C$.
    We set $\psi_P = \psi_{C,P}$.
    Let $H$ be a connected subgraph of $G\left[\bigcup\{Q \in \mathcal{P}, Q \geq_{\sigma} P\}\right]$ with $V(H) \cap P \neq \emptyset$.
    Then $H$ is a connected subgraph of $G_C\left[\bigcup\{Q \in \mathcal{P}_C\mid Q \geq_{\sigma_C} P\}\right]$.
    Hence by~\ref{item:thm:technical_Kt_minor_free_partition:parts_have_small_pcentered_colorings''} either $|\phi(V(H))|>p$ or $V(H)\cap P$ has a $(\phi\times\psi_{C,P})$-center $u$.
    In the former case, we are immediately satisfied, and in the latter case, $u$ is also a $(\phi\times\psi_P)$-center of $V(H)\cap P$.
    This yields~\ref{item:thm:technical_Kt_minor_free_partition:parts_have_small_pcentered_colorings} when $P \neq Z$.
    Finally, suppose that $P=Z$.
    Let $H$ be a connected subgraph of $G\left[\bigcup\{Q \in \mathcal{P}\mid Q \geq_{\sigma} Z\}\right] = G-R$.
    Now~\ref{item:thm:technical_Kt_minor_free_partition:parts_have_small_pcentered_colorings} follows directly from~\ref{items:pc-coloring-local-thm1-centers}.
    Summarizing, we constructed $\calP$, $\calR$, and $\sigma$ satisfying~\ref{item:thm:technical_Kt_minor_free_partition:Pi=Ri}-\ref{item:thm:technical_Kt_minor_free_partition:parts_have_small_pcentered_colorings}, as desired.
\end{proof}

\begin{proof}[Proof of \Cref{lemma:ISW-lifted}]
Let $t$ be a positive integer, and set $c_{\ref{lemma:ISW-lifted}}(t) = c_{\ref{lemma:Kt_free_graphs_have_centered_Helly_colorings}}(t) \cdot 2^{t-2}(t-1)$.
Let $G$ be a $K_t$-minor-free graph and $p$ be a positive integer.
By~\Cref{lemma:Kt_free_graphs_have_centered_Helly_colorings}, there exists a $(p,c_{\ref{lemma:Kt_free_graphs_have_centered_Helly_colorings}}(t))$-good coloring $\varphi$ of $G$ using at most $p+1$ colors.
Next, we apply~\Cref{thm:technical_Kt_minor_free_partition} with $c = c_{\ref{lemma:Kt_free_graphs_have_centered_Helly_colorings}}(t)$ and $r=0$.
We obtain a partition $\calP$ of $V(G)$, a tree decomposition $\calW$ of $G/\mathcal{P}$ of width at most $t-2$, and an elimination ordering $\sigma = (P_1, \dots, P_\ell)$ of $\calW$ such that for every $P \in \mathcal{P}$, there is a coloring $\psi_P\colon P \to [c_{\ref{lemma:Kt_free_graphs_have_centered_Helly_colorings}}(t) \cdot 2^{t-2}(t-1)]$ such that for every connected subgraph $H$ of $G\left[\bigcup\{Q \in \mathcal{P}\mid Q \geq_\sigma P\}\right]$ with $V(H) \cap P \neq \emptyset$, either $|\phi(V(H))| > p$, or $V(H) \cap P$ has a $\phi \times \psi_P$-center.
For every $P \in \calP$ and $u \in P$, let $\rho(u) = (\varphi(u), \psi_P(u))$.
Thus, $\rho$ is a coloring of $G$ using at most $c_{\ref{lemma:ISW-lifted}}(t)\cdot(p+1)$ colors and satisfying the assertion of~\Cref{lemma:ISW-lifted}.
\end{proof}

\section{Constructing good colorings}\label{sec:building_phi}

In this section, we prove \Cref{lemma:Kt_free_graphs_have_centered_Helly_colorings}, namely, we show that for every positive integer $t$, there exists a positive integer $c_{\ref{lemma:Kt_free_graphs_have_centered_Helly_colorings}}(t)$ such that for every $K_t$-minor-free graph $G$ and every positive integer $p$, $G$ admits a $(p,c_{\ref{lemma:Kt_free_graphs_have_centered_Helly_colorings}}(t))$-good coloring using $p+1$ colors.

At one point in the proof, we would like to assume that the tree decomposition that we consider is natural but still preserves certain properties.
For this reason, we introduce the following definitions and we prove~\Cref{lemma:making_a_td_natural}.

\begin{figure}[tp]
  \begin{center}
    \includegraphics{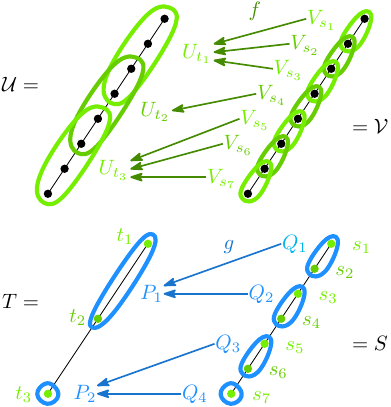}
  \end{center}
  \caption{
  In this example, $G$ is an $8$-vertex path.
  $\calP = \{P_1,P_2\}$ and $\calQ = \{Q_1,Q_2,Q_3,Q_4\}$.
  One can check that $(\calV,\calQ)$ refines $(\calU,\calP)$ 
  which is witnessed by $f$ and $g$.
  Note that $f$ is a function on $V(S)$ but for the readability reasons we depict it as it acted on bags of $\calV$. }
  \label{fig:refinement}
\end{figure}

Let $G$ be a graph. 
A pair $(\calU,\calP)$ is a \emph{normal pair} of 
$G$ if $\calU=\big(T,(U_t\mid t\in V(T))\big)$ is a tree decomposition of $G$, 
and $\calP$ is a collection of connected subgraphs 
of $T$ whose vertex sets partition $V(T)$.
Consider two normal pairs 
$(\calU,\calP)$ and $(\calV,\calQ)$ of $G$ with $\calU=\big(T,(U_t\mid t\in V(T))\big)$ and $\calV=\big(S,(V_s\mid s\in V(S))\big)$. 
We say that $(\calV,\calQ)$ \emph{refines}, or is a \emph{refinement} of $(\calU,\calP)$ if there exist $f\colon V(S) \rightarrow V(T)$ and $g\colon \calQ \rightarrow \calP$ such that
\begin{enumerate}[label=(r\arabic*)]
    \item for every $s \in V(S)$, 
        \[V_s \subset U_{f(s)};\] \label{item:making_a_td_natural:i}
    \item for every $Q \in \calQ$, 
        \[f(V(Q)) \subset V(g(Q)).\]\label{item:making_a_td_natural:ii}
\end{enumerate}
See an example in~\cref{fig:refinement}.

\begin{lemma}\label{lemma:making_a_td_natural}
    Let $G$ be a connected graph and $(\calU,\calP)$ be a normal pair of $G$. 
    There exists $(\calV,\calQ)$ a refinement of $(\calU,\calP)$ such that $\calV$ is natural.
\end{lemma}

\begin{proof}
For every positive integer $i$, for every tree decomposition $\mathcal{W}$ of $G$,
we denote by $n_i(\mathcal{W})$ the number of bags of $\mathcal{W}$ of size $i$.
Then let $n(\mathcal{W}) = (n_{|V(G)|}(\mathcal{W}), \dots, n_0(\mathcal{W}))$.

Since the refinement relation is reflexive, every normal pair of $G$ has a refinement.
Consider the refinement $(\calV,\calQ)$ of $(\calU,\calP)$ with $n(\calV)$ minimal in the lexicographic order.
We claim that $\calV$ is natural.

Let $\calV = \big(S, (V_s \mid s \in V(S))\big)$.
Suppose to the contrary that there exist $x,y \in V(S)$ such that $xy \in E(S)$
and $G\left[\bigcup_{s \in V(S_{x|y})} V_s\right]$ is not connected, and let $\mathcal{C}$ be the family of its components.

\begin{figure}[tp]
  \begin{center}
    \includegraphics{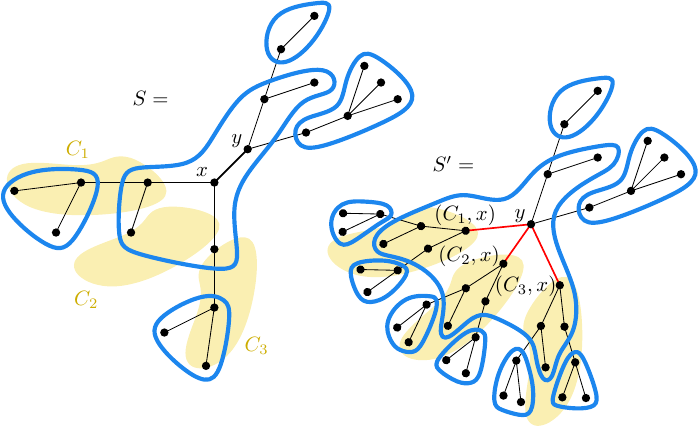}
  \end{center}
  \caption{
  Proof of~\Cref{lemma:making_a_td_natural}: construction of $S'$ from $S$ and of $\calQ'$ from $\calQ$. The blue partition of $S$ is $\calQ$ and the blue partition of $S'$ is $\calQ'$.
  }
  \label{fig:finidng-natural-td}
\end{figure}

The plan for obtaining a contradiction is finding a refinement $(\calV',\calQ')$ of $(\calV,\calQ)$ such that $n(\calV') < n(\calV)$ in the lexicographic order.
For every $C \in \mathcal{C}$, let $S_C$ be a copy of $S_{x|y}$
with vertex set $\{(C,s) \mid s \in V(S_{x|y})\}$ and edge set
$\{(C,s_1)(C,s_2) \mid s_1s_2 \in E(S_{x|y})\}$.
Let $S'$ be the tree obtained from the disjoint union of $S_{y|x}$ with all the $S_C$ over $C \in \mathcal{C}$ by adding the edges $(C,x)y$ for each $C \in \mathcal{C}$.
See~\cref{fig:finidng-natural-td}.
Let $f\colon V(S') \to V(S)$ be such that for every $z \in V(S')$,
\[
f(z) = \begin{cases} 
             s &\textrm{if $z=s$ for some $s\in V(S_{y|x})$,}\\
             s &\textrm{if $z=(C,s)$ for some $s\in V(S_{x|y})$ and $C\in\calC$.}
        \end{cases}
\]
Also, for every $z \in V(S')$, let
\[
V'_{z} = \begin{cases} 
V_{s}&\textrm{if $z=s$ for some $s\in V(S_{y|x})$,}\\
V_{s} \cap V(C)&\textrm{if $z=(C,s)$ for some $s\in V(S_{x|y})$ and $C\in\calC$.}
\end{cases}
\]
Note that by construction, for every $z\in V(S')$, $V'_z \subseteq V_{f(z)}$.
Moreover, $\calV'=\big(S',(V'_{z}\mid z\in V(S'))\big)$ is a tree decomposition of $G$.

To complete the construction, we specify a collection $\calQ'$ of subtrees of $S'$ whose vertex sets partition $V(S')$.
Simultaneously, we define $g\colon \calQ' \rightarrow \calQ$.
First, for every $Q \in \calQ$ with $Q \cap V(S_{y\mid x}) \neq \emptyset$, let 
    \[Q'=S'\left[(V(Q) \cap V(S_{y \mid x})) \cup \bigcup_{C \in \mathcal{C}} \{(C,z) \mid z \in V(Q) \cap V(S_{x \mid y})\}\right]\textrm{ and }g(Q')=Q.\]
Next, for every $Q \in \calQ$ with $Q \subseteq V(S_{x\mid y})$, for every $C \in \calC$, let
    \[Q_C = S'\left[\{(C,z) \mid z \in Q\}\right]\textrm{ and }g\big(Q_C\big)=Q.\]
Finally, let 
\begin{align*}
    \calQ' &= \left\{ Q' \mid Q \in \calQ \text{ with } V(Q) \cap V(S_{y\mid x}) \neq \emptyset  \right\} \cup \bigcup_{C \in \calC} \left\{ Q_C \mid V(Q) \subseteq V(S_{x\mid y})  \right\}.
\end{align*}
It follows that $\calQ'$ is a collection of subtrees of $S'$ whose sets of vertices partition $V(S')$ and 
and $f(V(Q')) = Q = g(Q)$ for every $Q' \in \calQ'$.
In other words,~\ref{item:making_a_td_natural:ii} holds.
Note that~\ref{item:making_a_td_natural:i} also holds, and thus, $(\calV',\calQ')$ refines $(\calV,\calQ)$.
Moreover, by transitivity of the refinement relation, $(\calV',\calQ')$ refines $(\calU,\calP)$.

To conclude the proof, we show that $n(\mathcal{V}')$ is less than $n(\mathcal{V})$ in the lexicographic order, which leads to a contradiction.
For every $s \in V(S_{x \mid y})$, we have $\sum_{C \in \calC} |V_{(C,s)}'| = |V_s|$. 
Thus, either $|V_{(C,s)}'| < |V_s|$ for every $C \in \calC$ or there is exactly one $C \in \calC$ with $|V_{(C,s)}'| = |V_s|$ and $|V_{(D,s)}'| = 0$ for all $D \in \calC\setminus\{C\}$.
Moreover, for every $s \in V(S_{y \mid x})$, we have $V_s = V_s'$.
Since $G$ is connected, $V_x$ intersects every component in $\mathcal{C}$.
Also, note that we supposed $|\calC| \geq 2$.
Consider the maximum integer $i_0$ such that there exists $s \in V(S_{x \mid y})$ 
with $|V_s|=i_0$ and $V_s$ intersects at least two components in $\mathcal{C}$.
From the construction, we obtain that $n_i(\mathcal{V}') = n_i(\mathcal{V})$ for every integer $i$ with $i > i_0$ and $n_{i_0}(\mathcal{V}')< n_{i_0}(\mathcal{V})$.
We conclude that $n(\mathcal{V}')$ is less than $n(\mathcal{V})$ in the lexicographic order.
\end{proof}

\begin{proof}[Proof of \Cref{lemma:Kt_free_graphs_have_centered_Helly_colorings}]
    Let $t$ be a positive integer. 
    Let $c_{\ref{theorem:Kt_free_product_structure_decomposition}}(t)$ be the constant from~\Cref{theorem:Kt_free_product_structure_decomposition}. 
    We set
    \[
    c_{\ref{lemma:Kt_free_graphs_have_centered_Helly_colorings}}(t) = c_{\ref{theorem:Kt_free_product_structure_decomposition}}(t) \left(12c_{\ref{theorem:Kt_free_product_structure_decomposition}}(t) + 10\right).
    \]
    
    Let $G$ be a $K_t$-minor-free graph, 
    and $p$ be a positive integer.
    By \Cref{theorem:Kt_free_product_structure_decomposition}, $G$ admits a layered RS-decomposition 
    $(T, \mathcal{W}, \mathcal{A}, \mathcal{D}, \mathcal{L})$ 
    of width at most $c_{\ref{theorem:Kt_free_product_structure_decomposition}}(t)$.
    Let $\mathcal{W} = (W_x \mid x \in V(T))$,
    $\mathcal{A} = (A_x \mid x \in V(T))$,
    $\mathcal{D} = \big(\big(T_x,(D_{x,z} \mid z \in V(T_x))\big) \mid x \in V(T)\big)$, and
    $\mathcal{L} = ((L_{x,i} \mid i \in \NN) \mid x \in V(T))$.
    We root $T$ in an arbitrary vertex $r \in V(T)$.

    We define a coloring $\phi \colon V(G) \rightarrow \{0,\dots,p\}$ as follows. For every $v\in V(G)$ let
    \[
    \phi(v) = \begin{cases}
    0&\textrm{if $v\in A_r$,}\\
    i\bmod (p+1)&\textrm{if $v\in W_r-A_r$ and $v \in L_{r,i}$,}\\
    0&\textrm{if $v\in A_x-W_{\parent(T,x)}$ for $x\in V(T)-\set{r}$,}\\
    i\bmod (p+1)&\textrm{if $v\in W_x-(A_x\cup W_{\parent(T,x)})$ and $v \in L_{x,i}$ for $x\in V(T)-\set{r}$.}\\
    \end{cases}
    \]

    It remains to  show that $\phi$ is a $(p,c_{\ref{lemma:Kt_free_graphs_have_centered_Helly_colorings}}(t))$-good coloring of $G$.
    Fix a subgraph $G_0$ of $G$. 
    
    Now, we reduce to the case where $G_0$ is connected. 
    Indeed, if $G_0$ is not connected, then consider the family $\mathcal{C}$ of all the components of $G_0$.
    Let $d$ be a nonnegative integer, and let $\mathcal{F}$ be a family of connected subgraphs of $G_0$ such that there are no $d+1$ disjoint members of $\mathcal{F}$.
    For every $C \in \mathcal{C}$, let $d_C$ be the smallest integer such that there are no $d_C+1$ disjoint members of $\mathcal{F}_C = \{F \in \mathcal{F} \mid F \subseteq C\}$.
    Clearly $d_C \leq d$. 
    Assuming we can prove the result when $G_0$ is connected,
    we apply it for $C$ and $\mathcal{F}_C$ and obtain a set $Z_C \subseteq V(C)$ and $\psi_C \colon Z_C \to [c \cdot d]$ such that \ref{centered-Helly-hitting}-\ref{centered-Helly-components} holds.
    Now take $Z = \bigcup_{C \in \mathcal{C}} Z_C$,
    and let $\psi$ be defined by $\psi(u) = \psi_C(u)$ for every $C \in \mathcal{C}$ and every $u \in Z_C$.
    Recall that each $F\in\calF$ is connected, so $F\subseteq C$ for some $C\in\calC$, and therefore $Z_C\cap V(F)\neq\emptyset$. 
    Thus,~\ref{centered-Helly-hitting} holds.
    Since every connected subgraph $H$ of $G_0$
    is a subgraph of $C$ for some $C \in \mathcal{C}$, 
    \ref{centered-Helly-coloring} holds. 
    Finally,~\ref{centered-Helly-components} holds because for every component $C'$ of $G-Z$, there exists $C \in \mathcal{C}$ such that $C' \subseteq C$.

    Therefore, from now on we assume that $G_0$ is connected.
    We start by building a normal pair $(\calU,\calP)$ where $\calU$ is a tree decomposition of $G_0$ obtained from ``gluing'' tree decompositions in $\calD$ along the edges of $T$.
    This construction is depicted in~\Cref{fig:new-td}.

\begin{figure}[tp]
  \begin{center}
    \includegraphics{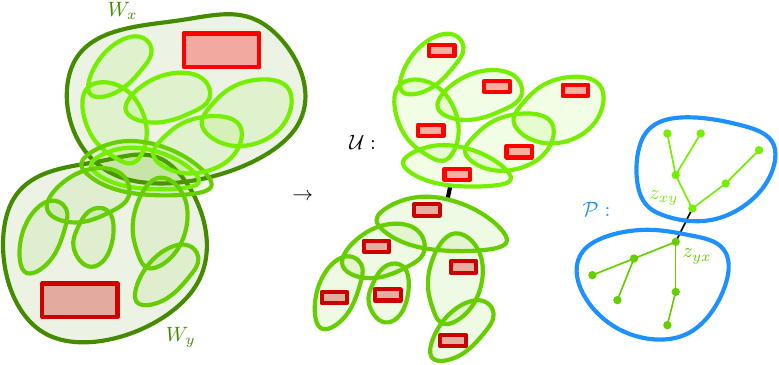}
  \end{center}
  \caption{
  Tree decompositions $\calD_x$ and $\calD_y$ are ``glued'' in a natural way. 
  Additionally, we add all corresponding apices (members of $A_x$ and $A_y$) to each bag of the corresponding part in the tree decomposition $\calU$.
  The partition $\calP$ indicates from which of $\{T_x \mid x \in V(T)\}$ a given vertex comes from.
  }
  \label{fig:new-td}
\end{figure}

    For all $x,y \in V(T)$ with $xy \in E(T)$, $W_x \cap W_y$ induces a clique in both $\torso_G(W_x)$ and $\torso_G(W_y)$, hence, there exist $z_{xy} \in V(T_x)$ and $z_{yx} \in V(T_y)$ such that
    $W_x \cap W_y = D_{x,z_{xy}} \cap D_{y,z_{yx}}$.
    Let $T_0$ be the tree defined by
    \begin{align*}
        V(T_0) &= \{(x,z) \mid x \in V(T), z \in V(T_x)\} \text{ and }\\
        E(T_0) &= \{(x,z_1)(x,z_2) \mid x \in V(T), z_1z_2 \in E(T_x)\} \cup \{(x, z_{xy}) (y,z_{yx}) \mid xy \in E(T)\}.
    \end{align*}
    In other words, $T_0$ is obtained from the disjoint union of $T_x$ over $x \in V(T)$ by adding the edges between $z_{xy}$ and $z_{yx}$ for adjacent vertices $x$ and $y$ in $T$.
    Next, let $U_{(x,z)} = (D_{x,z} \cup A_x) \cap V(G_0)$ for every $(x,z) \in V(T_0)$.
    We claim that $\mathcal{U} = \big(T_0, (U_{(x,z)} \mid (x,z) \in V(T_0))\big)$ is a tree decomposition of $G_0$.
    For every $v\in V(G_0)$, $T_0[\{(x,z)\in V(T_0)\mid v\in U_{(x,z)}\}]$ is isomorphic to the connected subgraph of $T_0$ formed by replacing every vertex $x$ of $T(v)=T[\{x\in V(T)\mid v\in W_x\}]$ with $T_x[\{z\in V(T_x)\mid v\in D_{x,z}\}]$ when $v\notin A_x$ or with $T_x$ when $v\in A_x$, and every edge $xy$ of $T(v)$ with $z_{xy}z_{yx}$.
    For every edge $uv\in E(G_0)$, there exists $x\in V(T)$ such that $u,v\in W_x$.
    If $u, v \in W_x \setminus A_x$, then there exists $z\in V(T_x)$ such that $u,v\in D_{x,z}$,
    and thus $u,v\in U_{(x,z)}$.
    If $u\in A_x$ and $v\in W_x \setminus A_x$, then for any $z\in V(T_x)$ with $v\in D_{x,z}$ we have $v\in U_{(x,z)}$ and $u\in A_x\cap V(G_0)\subseteq U_{(x,z)}$.
    If $u,v\in A_x$, then for every $z \in V(T_x)$, we have $u,v\in U_{(x,z)}$.
    Thus, indeed $\mathcal{U}$ is a tree decomposition of $G_0$.
    
    Observe that $\mathcal{P} = \{T_0[\{x\} \times V(T_x)] \mid x \in V(T)\}$ is a collection of subtrees of $T_0$ whose vertex sets partition $V(T_0)$. 
    Therefore, $\big(\mathcal{U},\mathcal{P}\big)$ is a normal pair of $G_0$.
    
    By \Cref{lemma:making_a_td_natural}, 
    there exists $\mathcal{V} = \big(S,(V_y \mid y \in V(S))\big)$
    and $\mathcal{Q}$ such that $\mathcal{V}$ is natural and $(\mathcal{V},\mathcal{Q})$ is a normal pair of $G$ refining $(\mathcal{U},\mathcal{P})$. 
    Among all such pairs $(\mathcal{V},\mathcal{Q})$, we take one with $|\mathcal{Q}|$ minimum.
    Let $f\colon V(S) \to V(T_0)$ and $g \colon \mathcal{Q} \to \mathcal{P}$ witness the refinement relation.
    For every $Q \in \mathcal{Q}$, let $x(Q) \in V(T)$ 
    be such that $g(Q) = T_0[\{x(Q)\} \times V(T_{x(Q)})]$.
    See~\Cref{fig:VQ}.

    In the next two claims, we show that the new tree decomposition $\calV$ in some sense preserves small adhesions (only ones coming from $\calW$) and small intersections with layers of the layerings in $\calL$.
    In the first claim, we exploit the minimality of $(\calV,\calQ)$. 
    \begin{claim}\label{claim:small_adhesions_between_Q}
        Let $y,y' \in V(S)$ and $Q,Q' \in \mathcal{Q}$ such that
        $y \in V(Q)$, $y' \in V(Q')$ and $Q \neq Q'$.
        Then
        \[
        |V_y \cap V_{y'}| \leq c_{\ref{theorem:Kt_free_product_structure_decomposition}}(t).
        \]
    \end{claim}

    \begin{proofclaim}
        By properties of tree decompositions,
        there is an edge $zz'$ in $S$ such that $V_y \cap V_{y'} \subseteq V_{z} \cap V_{z'}$.
        Therefore, without loss of generality, we assume that $yy'$ is an edge in $S$.
        We argue that $g(Q) \neq g(Q')$.
        Suppose to the contrary that $g(Q)=g(Q')$. 
        Since $yy'$ is an edge in $S$, the subgraphs $Q$ and $Q'$ are adjacent in $S$, and thus, $Q'' = S[V(Q) \cup V(Q')]$ is a connected subgraph of $S$.
        Observe that $(\mathcal{V}, \mathcal{Q} \setminus \{Q,Q'\} \cup \{Q''\})$ is a normal pair which refines $(\mathcal{U},\mathcal{P})$
        as witnessed by $f$ and $g'$ defined by $g'(P) = g(P)$ for every $P \in \mathcal{Q} \setminus \{Q,Q'\}$ and
        $g'(Q'') = g(Q) = g(Q')$.
        However, this contradicts the minimality of $|\mathcal{Q}|$.
        We obtain that indeed $g(Q) \neq g(Q')$, and so, $x(Q) \neq x(Q')$.
        In particular,
        \[
            |V_y \cap V_z|
            \leq |U_{(x(Q),f(y))} \cap U_{(x(Q'),f(y'))}| 
            \leq |W_{x(Q)} \cap W_{x(Q')}|  
            \leq c_{\ref{theorem:Kt_free_product_structure_decomposition}}(t),
        \]
        where the first inequality follows from~\ref{item:making_a_td_natural:i}, the second from the properties of tree decompositions, and the last from~\ref{LRS:adhesion}.
    \end{proofclaim}

\begin{figure}[tp]
  \begin{center}
    \includegraphics{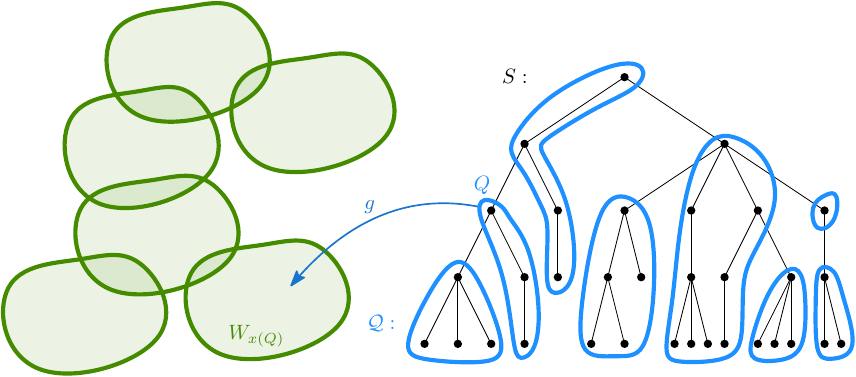}
  \end{center}
  \caption{
    Intuitively, the part of the tree decomposition $\calV$ corresponding to the vertices in $Q$ comes from $\calD_{x(Q)}$. 
  }
  \label{fig:VQ}
\end{figure}
    
    \begin{claim}\label{claim:Vy_intersection_Lxi_small}
    Let $x\in V(T)$, $y\in V(S)$, and $i\in \NN$. Then 
    \[
    |V_y\cap L_{x,i}|\leq  c_{\ref{theorem:Kt_free_product_structure_decomposition}}(t).
    \]
    \end{claim}
    \begin{proofclaim}
    Since $(\calV,\calQ)$ refines $(\calU,\calP)$ we have that 
    $V_y\subseteq U_{f(y)}$. Let $f(y)=(x',z)$ where $x'\in V(T)$ and $z\in V(T_{x'})$. Thus, 
    $U_{f(y)}=U_{(x',z)}=(D_{x',z}\cup A_{x'})\cap V(G_0)$.
    Recall that $A_x\cap L_{x,i}=\emptyset$. Thus in the case of $x=x'$,
    \[
    |V_y\cap L_{x,i}| 
    \leq |U_{(x,z)} \cap L_{x,i}| 
    \leq |(D_{x,z}\cup A_x)\cap L_{x,i}|
    \leq |D_{x,z}\cap L_{x,i}|
    \leq c_{\ref{theorem:Kt_free_product_structure_decomposition}}(t),
    \]
    where the last inequality follows by~\ref{LRS:ltw}.
    Finally, in the case of $x\neq x'$, we have 
    $U_{(x',z)}\subseteq W_{x'}$ and $L_{x,i}\subseteq W_{x}$, and so,
    \[
    |V_y\cap L_{x,i}| 
    \leq |U_{(x',z)} \cap L_{x,i}| 
    \leq |W_{x'}\cap W_x|
    \leq c_{\ref{theorem:Kt_free_product_structure_decomposition}}(t),
    \]
    where the last inequality follows from~\ref{LRS:adhesion}.
    \end{proofclaim}
    
    We root $S$ in an arbitrary vertex $\root(S)$.
    For every $X \subset V(S)$, we define
        \[\calQ(X) = \set{Q\in\calQ\mid\ V(Q) \cap X \neq \emptyset}.\]

    Let $d$ be a positive integer,
    and let $\mathcal{F}$ be a family of connected subgraphs of $G_0$
    such that there are no $d+1$ pairwise disjoint members of $\mathcal{F}$.
    In the remainder of the proof, we construct $Z$ and $\psi$ satisfying \ref{centered-Helly-hitting}-\ref{centered-Helly-components}.
    
    By \Cref{lemma:helly_property_tree_decomposition}, there exists $X_0 \subseteq V(S)$ of size at most $d$ such that $\bigcup_{y \in X_0} V_y$ intersects every member of $\mathcal{F}$.
    Let
        \[X_1 = \LCA(S,X_0).\] 
    By~\Cref{lemma:increase_X_in_a_tree},  
    \begin{equation}\label{eq:X1-and-Q1}
    |\calQ(X_1)|\leq |X_1| \leq 2d-1.
    \end{equation}

    The set $\bigcup_{y \in X_1} V_y$ is a hitting set of $\calF$, thus, it is a good candidate to be $Z$.
    However, we still need to define $\psi$ and the ultimate goal is to repeat the idea of coloring that we applied in the bounded layered treewidth case.
    Since $\calL_x$ is a layering of $\torso_G(W_x)-A_x$, we need to take into account the vertices in $A_x$ when building the final $Z$.
    Moreover, the coloring $\varphi$ may not be compatible with layerings on the adhesions of $\calW$, we also need to consider some of them.
    To this end, we will add some more vertices to $Z$, which we ultimately color injectively with a separate palette of colors.
    Because we work with $\calV$ instead of $\calU$ we need the following notions of projections. 
    See also~\Cref{fig:projections} for an illustration.

\begin{figure}[tp]
  \begin{center}
    \includegraphics{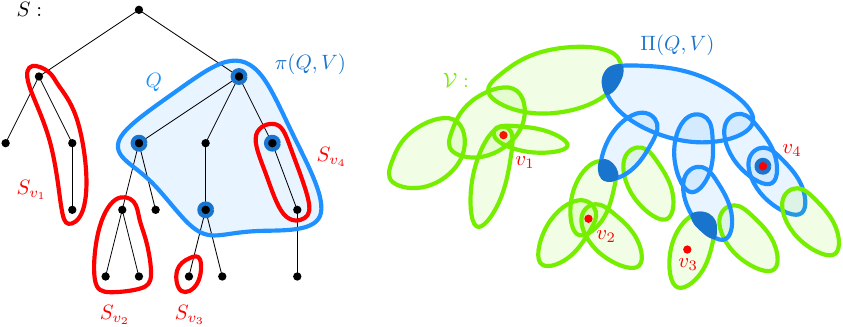}
  \end{center}
  \caption{
    Here, $V = \{v_1,v_2,v_3,v_4\}$.
    Note that $\bigcup_{y \in V(Q)} V_y$ is separated from $V$ by $\Pi(Q,V)$ in $G_0$ as we show later in~\Cref{claim:Pi-separates-stuff}. 
  }
  \label{fig:projections}
\end{figure}

    Let $Q \in \mathcal{Q}$
    and $v \in V(G_0)$.
    Let $S_v = S[\set{z\in V(S)\mid v\in V_z}]$.
    Thus, $Q$ and $S_v$ are two subtrees of $S$. 
    We define projections $\pi(Q,v)$ and $\Pi(Q,v)$ as follows.
    If $Q$ and $S_v$ share a vertex, then choose such a vertex $z$ arbitrarily and set $\pi(Q,v) = z$ and $\Pi(Q,v)=\set{v}$. 
    If $V(Q)$ and $V(S_v)$ are disjoint, then 
    consider the shortest path in $S$ between $V(Q)$ and $V(S_v)$. 
    Let $z$ be the endpoint of that path in $Q$ and 
    $y$ be the vertex adjacent to $z$ in that path. 
    Then set $\pi(Q,v) = z$ and $\Pi(Q,v) = V_y\cap V_{z}$.
    The definition can be naturally extended to subsets of $V(G_0)$.
    For each $V\subseteq V(G_0)$, let 
    \[\text{$\pi(Q,V)=\set{\pi(Q,v)\mid v\in V}$ \ and \ 
    $\Pi(Q,V)=\bigcup_{v \in V} \Pi(Q,v)$.}\]
    Note that by~\Cref{claim:small_adhesions_between_Q},
    \begin{equation}\label{eq:size_Pi}
        |\Pi(Q,V)| \leq c_{\ref{theorem:Kt_free_product_structure_decomposition}}(t) \cdot |V|.
    \end{equation}
    The key property of these objects is the following. 
    See~\Cref{fig:projections} again.
    \begin{claim}\label{claim:Pi-separates-stuff}
        Let $Q \in \calQ$ and $V \subset V(G_0)$.
        Every connected subgraph $H$ of $G_0 - \Pi(Q,V)$ that intersects $\bigcup_{y \in V(Q)} V_y$ is disjoint from $V$.
    \end{claim}
    \begin{proof}
        It suffices to observe that by the properties of tree decompositions and construction, $\Pi(Q,V)$ intersects every path between $\bigcup_{y \in V(Q)} V_y$ and $V$ in $G_0$.
        Thus, if a connected subgraph $H$ of $G_0$ has a vertex is both $\bigcup_{y \in V(Q)} V_y$ and $V$, then it has a vertex in $\Pi(Q,V)$.
    \end{proof}

    Let
    \[
    Y_Q =
    \begin{cases}
        \pi(Q, A_{x(Q)}) &\textrm{if $x(Q)=r$,}\\
        \pi(Q, A_{x(Q)} \cup (W_{x(Q)} \cap W_{\parent(T,x(Q))})) & \textrm{if $x(Q) \neq r$.}
    \end{cases}
    \]
    Recall that $|A_{x(Q)}|\leq c_{\ref{theorem:Kt_free_product_structure_decomposition}}(t)$ by~\ref{LRS:apices},
    and if $x(Q)\neq r$, then $|W_{x(Q)} \cap W_{\parent(T,x(Q))}|\leq c_{\ref{theorem:Kt_free_product_structure_decomposition}}(t)$ by~\ref{LRS:adhesion}.
    Thus, we have 
    \begin{equation}
    \begin{aligned}
    |Y_Q|
    & \leq |A_{x(Q)}| 
    &&\leq c_{\ref{theorem:Kt_free_product_structure_decomposition}}(t) && \textrm{if $x(Q)=r$,}\\
    |Y_Q|
    & 
    \leq |A_{x(Q)}|+|W_{x(Q)}\cap W_{\parent(T,x(Q))}| 
    &&\leq 2c_{\ref{theorem:Kt_free_product_structure_decomposition}}(t) && \textrm{if $x(Q)\neq r$.}
    \end{aligned}
    \label{eq:YQ}
    \end{equation}
    Moreover, by definition $Y_Q \subset V(Q)$.

    Let
    \[
        X_2 = X_1 \cup \{\root(Q) \mid Q \in \mathcal{Q}(X_1)\} \cup \bigcup \{ Y_Q \mid Q \in \mathcal{Q}(X_1)\},
    \]
    and
    \[
        X_3 = \LCA(S,X_2).
    \]
    Observe that
    \begin{equation}
    \begin{aligned}
    |X_3|
    &\leq 2|X_2| - 1 
    &&\textrm{by~\Cref{lemma:increase_X_in_a_tree}}\\
    &\leq 2(|X_1| + |\calQ(X_1)| + |\calQ(X_1)|\cdot 2c_{\ref{theorem:Kt_free_product_structure_decomposition}}(t))-1
    &&\textrm{by~\eqref{eq:YQ}}\\
    &\leq 2((2d-1) + (2d-1) + (2d-1)\cdot 2c_{\ref{theorem:Kt_free_product_structure_decomposition}}(t))-1 &&\textrm{by~\eqref{eq:X1-and-Q1}}\\
    &\leq (8c_{\ref{theorem:Kt_free_product_structure_decomposition}}(t)+8)d.
    \end{aligned}\label{eq:size-X_3}
    \end{equation}

    Recall that by~\Cref{lemma:LCA-stuff}, for all $X, Y \subseteq V(S)$ with $X\subseteq Y$ and $\LCA(S,X) = X$, if $\calQ(X) = \calQ(Y)$, then $\calQ(Y)=\calQ(\LCA(S,Y))$.
    Since by construction $\calQ(X_1) = \calQ(X_2)$, we can apply the above with $X = X_1$ and $Y = X_2$ to obtain that
    \[
        \mathcal{Q}(X_1) = \calQ(X_2) = \mathcal{Q}(X_3).
    \]

    Finally, the set $Z$ that witnesses the assertion of the claim is given by
    \[
    Z = \bigcup_{y \in X_3} V_y.
    \]
    
    Next, consider $Q \in \mathcal{Q}(X_1)$,
    and let 
    \[
    B_Q =
    \begin{cases}
        \Pi(Q, A_{x(Q)}) &\textrm{if $x(Q)=r$,} \\
        \Pi(Q, A_{x(Q)} \cup (W_{x(Q)} \cap W_{\parent(T,x(Q))})) & \textrm{if $x(Q) \neq r$.}
    \end{cases}
    \]
    Since for every $v \in V(G_0)$, $\Pi(Q,v) \subseteq V_{\pi(Q,v)}$ and $Y_Q \subseteq X_2 \subseteq X_3$, we have 
    \[
    B_Q \subseteq \bigcup_{y \in Y_Q} V_y \subseteq \bigcup_{y \in X_2} V_y \subseteq Z.
    \]
    Recall that $|A_{x(Q)}|\leq c_{\ref{theorem:Kt_free_product_structure_decomposition}}(t)$ by~\ref{LRS:apices},
    and if $x(Q)\neq r$, then $|W_{x(Q)} \cap W_{\parent(T,x(Q))}|\leq c_{\ref{theorem:Kt_free_product_structure_decomposition}}(t)$ by~\ref{LRS:adhesion}.
    Thus, by \eqref{eq:size_Pi}, we have
    \begin{equation}
    \begin{aligned}
    |B_Q|
    &\leq |A_{x(Q)}| \cdot c_{\ref{theorem:Kt_free_product_structure_decomposition}}(t)
    &&\leq (c_{\ref{theorem:Kt_free_product_structure_decomposition}}(t))^2
    && \textrm{if $x(Q)=r$,}\\
    |B_Q|
    &\leq \left|A_{x(Q)}\cup(W_{x(Q)}\cap W_{\parent(T,x(Q))}))\right| \cdot c_{\ref{theorem:Kt_free_product_structure_decomposition}}(t) 
    &&\leq 2(c_{\ref{theorem:Kt_free_product_structure_decomposition}}(t))^2
    &&\textrm{if $x(Q) \neq r$.}
    \end{aligned}
    \label{eq:A'Q}
    \end{equation}
    
    Let 
    \[
        B = \left(\bigcup_{Q \in \mathcal{Q}(X_1), \root(S) \not\in V(Q)} V_{\root(Q)} \cap V_{\parent(S,\root(Q))}\right) \cup \left(\bigcup_{Q \in \mathcal{Q}(X_1)} B_Q\right).
    \]
    Recall that for every $Q \in \mathcal{Q}(X_1)$ with $\root(S) \not\in V(Q)$,
    $\root(Q) \in X_2$ and so $V_{\root(Q)} \cap V_{\parent(S,\root(Q))} \subseteq V_{\root(Q)} \subseteq Z$. 
    Moreover, for every $Q \in \mathcal{Q}(X_1)$, $B_Q \subseteq Z$.
    Therefore,
    \[
    B\subseteq Z. 
    \]
    Observe also that
    \begin{equation}
    \begin{aligned}
    |B| 
    &\leq |\mathcal{Q}(X_1)|\cdot c_{\ref{theorem:Kt_free_product_structure_decomposition}}(t)  + |\mathcal{Q}(X_1)|\cdot 2(c_{\ref{theorem:Kt_free_product_structure_decomposition}}(t))^2 && \textrm{by \Cref{claim:small_adhesions_between_Q} and \eqref{eq:A'Q}} \\
    &\leq (2d-1)\cdot c_{\ref{theorem:Kt_free_product_structure_decomposition}}(t) 
    + (2d-1)\cdot 2(c_{\ref{theorem:Kt_free_product_structure_decomposition}}(t))^2
    &&\textrm{by \eqref{eq:X1-and-Q1}}\\
    & \leq c_{\ref{theorem:Kt_free_product_structure_decomposition}}(t)(4c_{\ref{theorem:Kt_free_product_structure_decomposition}}(t)+2)d.
    \end{aligned}
    \label{eq:size-of-B}
    \end{equation}
    
    \begin{claim}\label{claim:B_separates_the_Qs}
        For all distinct $Q_1,Q_2 \in \mathcal{Q}(X_1)$,
        there is no path between $\bigcup_{y \in V(Q_1)} V_y \setminus B$ and
        $\bigcup_{y \in V(Q_2)} V_y \setminus B$ in $G_0 - B$.
        In particular, $\bigcup_{y \in V(Q_1)} V_y \setminus B$ and $\bigcup_{y \in V(Q_2)} V_y \setminus B$ are disjoint.
    \end{claim}

    \begin{proofclaim}
        Let $Q_1,Q_2$ be distinct members of $\mathcal{Q}(X_1)$.
        Recall that $Q_1,Q_2$ are disjoint subtrees of $S$.
        Note that there exists $i \in \{1,2\}$
        such that $\root(Q_i) \neq \root(S)$ and every path in $S$ between $Q_1$ and $Q_2$
        goes through the edge $\root(Q_i)\parent(S,\root(Q_i))$ of $S$.
        Hence by properties of tree decompositions, 
        $V_{\root(Q_i)} \cap V_{\parent(S,\root(Q_i))}$ intersects 
        every path between $\bigcup_{y \in V(Q_1)} V_y$ and $\bigcup_{y \in V(Q_2)} V_y$ in $G_0$.
        Since $V_{\root(Q_i)} \cap V_{\parent(S,\root(Q_i))} \subseteq B$, this proves the claim.
    \end{proofclaim}

    As a consequence of \Cref{claim:B_separates_the_Qs},
    $\big\{\big(\bigcup_{y \in V(Q)} V_y\big) \cap (Z \setminus B) \mid Q \in \mathcal{Q}(X_1)\big\}$ is a family of pairwise disjoint sets covering $Z \setminus B$ (because $\mathcal{Q}(X_1) = \mathcal{Q}(X_3)$).
    We want to refine this family 
    with the layerings of the torsos.
    To this end, we need the following fact where we substantially use the fact that $(\calV,\calQ)$ refines $(\calU,\calP)$.
    \begin{claim}\label{claim:V(Q)-stays-in-one-W}
        For every $Q \in \calQ(X_1)$ and $V \subset V(G_0)$ such that $\Pi(Q,V) \subset B_Q$, we have
            \[\bigcup_{y \in V(Q)} V_y \setminus B \subseteq W_{x(Q)} \setminus V.\]
    \end{claim}
    \begin{proofclaim}
        Let $Q \in \calQ(X_1)$ and let $V \subset V(G_0)$ be such that $\Pi(Q,V) \subset B_Q$.
        Since $\Pi(Q,V) \subset B_Q \subset B$, then
            \[\bigcup_{y \in V(Q)} V_y \setminus B \subset \bigcup_{y \in V(Q)} V_y \setminus B_Q \subset \bigcup_{y \in V(Q)} V_y \setminus \Pi(Q,V).\]
        By~\cref{claim:Pi-separates-stuff}, $\bigcup_{y \in V(Q)} V_y \setminus \Pi(Q,V)$ is disjoint from $V$.
        Thus, to conclude the claim, it suffices to show that $\bigcup_{y \in V(Q)} V_y \subset W_{x(Q)}$.
        Recall that $g(Q) = T_0[\{x(Q)\}\times T_{x(Q)}]$.
        We have
            \[\bigcup_{y \in V(Q)} V_y \subset \bigcup_{y \in V(Q)} U_{f(y)} \subset \bigcup_{z \in V(g(Q))} U_z = \bigcup_{z \in \{x(Q)\} \times T_{x(Q)}} U_z \subset W_{x(Q)}\]
        where the first inclusion follows from~\ref{item:making_a_td_natural:i}, the second from~\ref{item:making_a_td_natural:ii}, and the last one from the construction of $\calU$.
    \end{proofclaim}

    Since  $\Pi(Q,A_{x(Q)}) \subset B_Q$ for every $Q \in \calQ(X_1)$, by~\cref{claim:V(Q)-stays-in-one-W},
        \[\textstyle \big(\bigcup_{y \in V(Q)} V_y - B\big) \cap A_{x(Q)} = \emptyset. \]
    Recall that for every $Q \in \calQ(X_1)$, $(L_{x(Q),i} \mid  i \in \NN)$ is a layering of
    $\torso_G(W_{x(Q)}) - A_{x(Q)}$.
    It follows that the family
    \[\textstyle
    \calB = \left\{\big(\bigcup_{y \in V(Q)} V_y\big) \cap L_{x(Q),i} \cap (Z \setminus B) \mid Q \in \mathcal{Q}(X_1), i \in \NN \right\}
    \]
    is a family of pairwise disjoint sets covering $Z \setminus B$.
    Moreover, the members of $\calB$ are of reasonable size, namely, for every $Q \in \calQ(X_1)$ and $i \in \NN$,
    \begin{equation} 
    \begin{aligned}
    \textstyle
        \left|\big(\bigcup_{y \in V(Q)} V_y\big) \cap L_{x(Q),i} \cap (Z \setminus B)\right| &\leq |L_{x(Q),i} \cap Z| \leq \sum_{y \in X_3} |L_{x(Q),i} \cap V_y|\\ 
         &\leq \sum_{y \in X_3} |L_{x(Q),i} \cap U_{f(y)}| && \text{by~\ref{item:making_a_td_natural:i}}\\
        &\leq c_{\ref{theorem:Kt_free_product_structure_decomposition}}(t) \cdot |X_3| && \text{by~\ref{LRS:ltw}}\\
        &\leq c_{\ref{theorem:Kt_free_product_structure_decomposition}}(t) \cdot \big(8c_{\ref{theorem:Kt_free_product_structure_decomposition}}(t)+8\big)d && \text{by~\eqref{eq:size-X_3}.}
    \end{aligned}
    \label{eq:bound_on_layer_cap_Z}
    \end{equation}

\begin{figure}[tp]
  \begin{center}
    \includegraphics{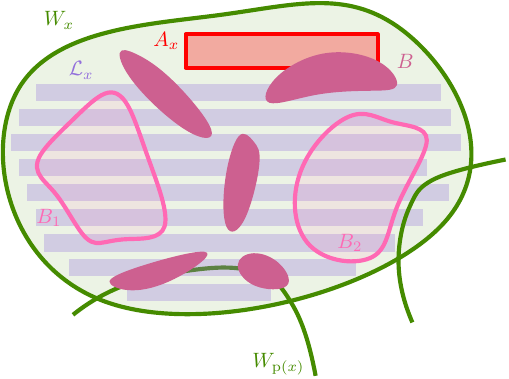}
  \end{center}
  \caption{
    The set $Z$ is the union of $B$ and the members of $\calB$.
    We color the vertices in $B$ injectively.
    The role of the set $B$ is to separate members of $\calB$ so that if $B$ is intersected by a connected subgraph $H$ of $G_0$, then we immediately get a center.
    In the figure, $B_1 = \bigcup_{y \in V(Q_1)} V_y - B$ and $B_2 = \bigcup_{y \in V(Q_2)} V_y - B$ for some $Q_1,Q_2 \in \calQ(X_3)$.
    If a connected subgraph $H$ of $G_0$, intersects both $B_1$ and $B_2$, then it has to intersect $B$.
  }
  \label{fig:constructing-psi}
\end{figure}

    We define a coloring $\psi$ of $Z$ using at most $c_{\ref{lemma:Kt_free_graphs_have_centered_Helly_colorings}}(t)$ colors.
    See also~\Cref{fig:constructing-psi}.
    First, we color $B$ injectively, and then we color each member of $\calB$ also injectively avoiding colors in $\psi(B)$.
    In the first step we used at most $c_{\ref{theorem:Kt_free_product_structure_decomposition}}(t)(4c_{\ref{theorem:Kt_free_product_structure_decomposition}}(t)+2)d$ colors by~\eqref{eq:size-of-B}, and in the second step, we used at most $c_{\ref{theorem:Kt_free_product_structure_decomposition}}(t) \cdot \big(8c_{\ref{theorem:Kt_free_product_structure_decomposition}}(t)+8\big)d$ colors by~\eqref{eq:bound_on_layer_cap_Z}, thus, $\psi$ is well-defined.

    We now show that $Z$ and $\psi$ satisfy \ref{centered-Helly-hitting}-\ref{centered-Helly-components}.
    
    Recall that $X_0$ was chosen so that $X_0 \subseteq X_3$ and $\bigcup_{y \in X_0} V_y$ intersects every member of $\mathcal{F}$.
    Therefore,~\ref{centered-Helly-hitting} holds.
    Item~\ref{centered-Helly-components} holds
    by \Cref{lemma:increase_X_to_have_small_interfaces} since $Z = \bigcup_{y \in \LCA(S, X_2)} V_y$, and $\mathcal{V} = \big(S,(V_y \mid y \in V(S))\big)$ is a natural tree decomposition of $G_0$.

    It remains to prove~\ref{centered-Helly-coloring}.
    Consider a connected subgraph $H$ of $G_0$ such that $V(H) \cap Z \neq \emptyset$.
    If $V(H) \cap B \neq \emptyset$,
    then since every vertex of $B$ has a unique color, any vertex in $V(H) \cap B$ is a $\psi$-center of $V(H) \cap Z$, and so, a $\phi \times \psi$-center of $V(H) \cap Z$.
    Thus, we assume that $V(H) \cap B = \emptyset$.
    
    Since $V(H) \cap B = \emptyset$, it follows from~\Cref{claim:B_separates_the_Qs}, that $V(H) \cap Z \subseteq \bigcup_{z \in V(Q)} V_z$ for some $Q \in \mathcal{Q}(X_1)$.

    Additionally, by \Cref{claim:V(Q)-stays-in-one-W},
    \begin{align*}
        V(H) \cap Z &\subseteq W_x \setminus A_x && \textrm{if $x = r$,}\\
        V(H) \cap Z &\subseteq W_x \setminus (A_x \cup W_{\parent(T,x)}) && \textrm{if $x \neq r$.}
    \end{align*}
    Recall that $(L_{x,i} \mid i \in \mathbb{N})$ is a layering of $\torso_G(W_x)-A_x$.
    Moreover, by definition of $\phi$,
    for every $u \in V(H) \cap Z$, $\phi(u) = i \bmod (p+1)$ where $i \in \NN$ is such that $u \in L_{x,i}$.

    Consider $u \in V(H) \cap Z$.
    Let $i\in \NN$ be such that $u \in L_{x,i}$.
    If $u$ is a $(\phi \times \psi)$-center in $V(H) \cap Z$, then we are done.
    Assume now that there exists $u' \in V(H) \cap Z$ distinct from $u$
    such that $\phi(u)=\phi(u')$ and $\psi(u) = \psi(u')$.
    Let $j \in \NN$ be such that $u' \in L_{x,j}$.
    Without loss of generality, assume that $i \leq j$.
    By the definition of $\psi$, $j \neq i$,
    and by the definition of $\phi$, $|j-i| > p$.
    By \Cref{lemma:projections_on_a_torso_stay_connected},
    $V(H) \cap W_x$ induces 
    a connected subgraph of $\torso_G(W_x)-A_x$.
    Since $(L_{x,k} \mid k \in \NN)$ is a layering of $\torso_G(W_x)-A_x$,
    it follows that $V(H)$ intersects $L_{x,k}$ for every $k \in \{i,\dots, j-1\}$. 
    We deduce that $|\phi(V(H))|>p$.
    Therefore, we obtain~\ref{centered-Helly-coloring}, which ends the proof.
\end{proof}

\bibliographystyle{plain}
\bibliography{biblio}

\begin{thebibliography}{10}

\bibitem{And86}
Thomas Andreae.
\newblock On a pursuit game played on graphs for which a minor is excluded.
\newblock {\em Journal of Combinatorial Theory, Series B}, 41(1):37--47, 1986.

\bibitem{DHHJLMMRW24}
Vida Dujmovic, Robert Hickingbotham, Jędrzej Hodor, Gwena{\"{e}}l Joret, Hoang La, Piotr Micek, Pat Morin, Cl{\'{e}}ment Rambaud, and David~R. Wood.
\newblock The grid-minor theorem revisited.
\newblock In David~P. Woodruff, editor, {\em Proceedings of the 2024 {ACM-SIAM} Symposium on Discrete Algorithms, {SODA} 2024, Alexandria, VA, USA, January 7-10, 2024}, pages 1241--1245. {SIAM}, 2024.
\newblock \href{https://arxiv.org/abs/2307.02816}{arXiv:2307.02816}.

\bibitem{DJMMUW20}
Vida Dujmovi\'{c}, Gwena\"{e}l Joret, Piotr Micek, Pat Morin, Torsten Ueckerdt, and David~R. Wood.
\newblock Planar graphs have bounded queue-number.
\newblock {\em Journal of the ACM}, 67(4), 2020.
\newblock \href{http://arxiv.org/abs/1904.04791}{arXiv:1904.04791}.

\bibitem{Dujmovi2017}
Vida Dujmović, Pat Morin, and David~R. Wood.
\newblock Layered separators in minor-closed graph classes with applications.
\newblock {\em Journal of Combinatorial Theory, Series B}, 127:111–147, 2017.
\newblock \href{https://arxiv.org/abs/1306.1595}{arXiv:1306.1595}.

\bibitem{Dbski2021}
Michał Dębski, Piotr Micek, Felix Schr\"{o}der, and Stefan Felsner.
\newblock Improved bounds for centered colorings.
\newblock {\em Advances in Combinatorics}, 2021.
\newblock \href{https://arxiv.org/abs/1907.04586}{arXiv:1907.04586}.

\bibitem{Eppstein1999}
David Eppstein.
\newblock Subgraph isomorphism in planar graphs and related problems.
\newblock {\em Journal of Graph Algorithms and Applications}, 3(3):1–27, 1999.

\bibitem{Grohe2015}
Martin Grohe and Dániel Marx.
\newblock Structure theorem and isomorphism test for graphs with excluded topological subgraphs.
\newblock {\em SIAM Journal on Computing}, 44(1):114–159, 2015.
\newblock \href{https://arxiv.org/abs/1111.1109}{arXiv:1111.1109}.

\bibitem{vdHetal17}
Jan {\noopsort{Heuvel}}{van den Heuvel}, Patrice {Ossona de Mendez}, Daniel Quiroz, Roman Rabinovich, and Sebastian Siebertz.
\newblock On the generalised colouring numbers of graphs that exclude a fixed minor.
\newblock {\em European Journal of Combinatorics}, 66:129--144, 2017.
\newblock \href{https://arxiv.org/abs/1602.09052}{arXiv:1602.09052}.

\bibitem{ISW22}
Freddie Illingworth, Alex Scott, and David~R. Wood.
\newblock Product structure of graphs with an excluded minor.
\newblock {\em arXiv preprint}, 2021.
\newblock \href{https://arxiv.org/abs/2104.06627}{arXiv:2104.06627}.

\bibitem{sparsity}
Jaroslav Ne\v{s}et\v{r}il and Patrice {Ossona de Mendez}.
\newblock {\em Sparsity --- {G}raphs, {S}tructures, and {A}lgorithms}, volume~28 of {\em Algorithms and {C}ombinatorics}.
\newblock Springer, 2012.

\bibitem{NOdMW12}
Jaroslav Nešetřil, Patrice~Ossona de~Mendez, and David~R. Wood.
\newblock Characterisations and examples of graph classes with bounded expansion.
\newblock {\em European Journal of Combinatorics}, 33(3):350 -- 373, 2012.

\bibitem{Nesetril2008}
Jaroslav Nešetřil and Patrice Ossona~de Mendez.
\newblock {Grad and classes with bounded expansion I. Decompositions}.
\newblock {\em European Journal of Combinatorics}, 29(3):760–776, April 2008.

\bibitem{notes}
Marcin Pilipczuk, Micha\l{} Pilipczuk, and Sebastian Siebertz.
\newblock Lecture notes for the course ``{S}parsity'' given at {F}aculty of {M}athematics, {I}nformatics, and {M}echanics of the {U}niversity of {W}arsaw, Winter semesters 2017/18 and 2019/20.
\newblock \url{https://www.mimuw.edu.pl/~mp248287/sparsity2}.

\bibitem{PS19}
Michal Pilipczuk and Sebastian Siebertz.
\newblock Polynomial bounds for centered colorings on proper minor-closed graph classes.
\newblock {\em {Journal of Combinatorial Theory, Series B}}, 151:111--147, 2021.
\newblock \href{https://arxiv.org/abs/1807.03683}{arXiv:1807.03683}.

\bibitem{GM5}
Neil Robertson and Paul Seymour.
\newblock Graph minors. {V}. {Excluding a planar graph}.
\newblock {\em Journal of Combinatorial Theory, Series B}, 41(1):92--114, 1986.

\bibitem{GM16}
Neil Robertson and Paul Seymour.
\newblock {Graph ninors. XVI. Excluding a non-planar graph}.
\newblock {\em Journal of Combinatorial Theory, Series B}, 89(1):43--76, 2003.

\end{thebibliography}

\end{document}